\documentstyle[11pt]{article}

\newcommand{\integ}[2]{\displaystyle \int_{#1}^{#2}}

\newcommand{\somme}[2]{\displaystyle \sum_{#1}^{#2}}
\newtheorem{theorem}{Theorem}[section]
\newtheorem{proposition}{Proposition}[section]
\newtheorem{lemma}{Lemma}[section]
\newtheorem{remark}{Remark}[section]
\newtheorem{corollary}{Corollary}[section]

\def \R{I\!\!R}

\def \T{\mathcal{T}}
\def \H{\mathcal{H}}
\def \L{\mathcal{L}}

\def \bf{\textbf}
\def \it{\textit}

\def \no {\noindent}

\title{BSDEs with two RCLL Reflecting Obstacles driven by a Brownian Motion and Poisson Measure and related
Mixed Zero-Sum Games}
\author {S.Hamad\`ene\, and\, H.Wang\\
Laboratoire de Statistique et Processus,\\ Universit\'e du Maine,
72085 Le Mans Cedex 9, France\\ {\small \it{e-mails}:
hamadene@univ-lemans.fr }\& {\small hao.wang.etu@univ-lemans.fr }}
\begin{document}

\date{\today}
\maketitle
\begin{abstract}
In this paper we study Backward Stochastic Differential Equations
with two reflecting right continuous with left limits obstacles (or
barriers) when the noise is given by Brownian motion and a Poisson
random measure mutually independent. The jumps of the obstacle
processes could be either predictable or inaccessible. We show
existence and uniqueness of the solution when the barriers are
completely separated and the generator uniformly Lipschitz. We do
not assume the existence of a difference of supermartingales between
the obstacles. As an application, we show that the related mixed
zero-sum differential-integral game problem has a value. 
\end{abstract}
\no \bf{AMS Classification subjects}: 91A15, 91B74, 60G40, 91A60

\no \textbf{Key Words:} Backward stochastic differential equation ;
penalization ; Mokobodski's hypothesis ; Snell envelope ; zero-sum
mixed differential-integral game.

\section{ Introduction} In this paper we are concerned with the
problem of existence and uniqueness of a solution for the backward
stochastic differential equations (BSDEs for short) driven by a
Brownian motion and an independent Poisson measure with two
reflecting obstacles (or barriers) which are right continuous with
left limits ($rcll$ for short) processes. Roughly speaking we look
for a quintuple of adapted processes $(Y,Z,V,K^\pm)$ such that:
\begin{equation}\label{bsdeintro}
      \left\{\begin{array}{l}
      (i)\,Y_t=\xi+\int_t^{T}f(s,Y_s,Z_s,V_s)ds+(K^+_T-K^+_t)-
      (K^-_T-K^-_t)-\int_t^{T}Z_sdB_s -
   \int_t^{T}\int_EV_s(e)\widetilde{\mu}(ds,de)\\
   (ii)\,L\leq Y\leq U \mbox{ and if }K^{c,\pm} \mbox{ is the
   continuous part of }K^\pm \mbox{ then
   }\\\qquad\qquad\int_0^T(Y_t-L_t)dK^{c,+}_t=\int_0^T(U_t-Y_t)dK^{c,-}_t=0\\
   (iii)\,\mbox{if }K^{d,\pm} \mbox{ is the purely discontinuous part
   of }K^\pm \mbox{ then }K^{d,\pm} \mbox{ is predictable and }\\
   \qquad K^{d,-}_t=\sum\limits_{0<s\leq t}(Y_s-U_{s-})^+1_{[\Delta
                                  U_s>0]} \mbox{ and }K^{d,+}_t=\sum\limits_{0<s\leq
                                  t}(Y_s-L_{s-})^-1_{[\Delta L_s<0]}\,\,
              \end{array}
      \right.
\end{equation}
where $B$ is a Brownian motion, $\tilde \mu$ is a compensated
Poisson random measure and $f(t,\omega,y,z,v)$, $\xi$, $L$ and $U$
are given ($B$ and $\tilde \mu$ are independent).

 In the framework of a Brownian filtration, the notion of BSDEs with one reflecting
 obstacle is introduced by El-Karoui et al. \cite{[EKal]}. Those equations
 have been well considered during the last ten
years since they have found a wide range of applications especially
in finance, stochastic control/games, partial differential
equations,.... Later Cvitanic $\&$ Karatzas generalized in
\cite{[CK]} the setting of \cite{[EKal]} where they introduced BSDEs
with two reflecting barriers. Since then there were several articles
on this latter types of BSDEs (see $e.g.$ \cite{[BHM],
 hh2, hibtis, [HL1], [HLM], [LSM], tmx, [PX]} and the
references therein), usually in connection with various
applications. Nevertheless during several years, the existence of a
solution of two barrier reflected BSDEs is obtained under one of the
two following hypotheses: either one of the obstacles is "almost" a
semimartingale (see $e.g.$ \cite{[CK],[HLM]}) or the so-called
Mokobodski's condition (see (\ref{mokocond}) for its definition)
\cite{[CK],[HL1],[LSM],tmx,[PX]} holds. Obviously the first
assumption is somehow restrictive as for the second one it is quite
difficult to check in practice. Those conditions have been removed
in \cite{hh2} where the authors showed that if the barriers are
continuous and completely separated, $i.e.$ $\forall t\leq
T,\,\,L_t<U_t$, then the two barrier reflected BSDE has a solution.
Later the case of discontinuous barriers has been also studied in
Hamad\`ene et al. \cite{hho} where they actually show the existence
of a solution when the obstacles and their left limits are
completely separated.

In this work, we focus on BSDEs with two reflecting barriers when,
on the one hand, the filtration is generated by a Brownian motion
and an independent Poisson random measure and, on the other hand,
the barriers are $rcll$ processes whose jumps are arbitrary, they
can be either predictable or inaccessible. We show that when the
generator of the BSDE is Lipschitz, the obstacle processes and their
left limits are completely separated then the BSDE (\ref{bsdeintro})
has a unique solution. Therefore our work is an extension of the one
by Hamad\`ene $\&$ Hassani \cite{hh3} where they deal with the same
framework of BSDEs except that the obstacle processes are not
allowed to have predictable jumps. This work generalizes also the
paper in \cite{hho} where the two barrier reflecting BSDE they
consider is driven only by a Brownian motion. The main difficulty of
our problem lies in the fact that the jumps of the obstacles can be
predictable or inaccessible, therefore the component $Y$ of the
solution has also both types of jumps. This is the basic difference
of our work related to \cite{hh3} (resp. \cite{hho}) where $Y$ has
only inaccessible (resp. predictable) jumps.

It is well known that double barrier reflected BSDEs are connected
with mixed zero-sum games (see $e.g.$ \cite{[H2],[HL1]}). Therefore
as an application of our result obtained in the first part of the
paper, in the second part we deal with zero-sum mixed stochastic
differential-integral games which we describe briefly. Assume we
have a system on which intervene two agents (or players) $c_1$ and
$c_2$. This system could be a stock in the market and then $c_1$,
$c_2$ are two traders whose advantages are antagonistic. The
intervention of the agents have two forms, control and stopping. The
dynamics of the system when controlled is given by:
$$\begin{array}{ll}
 x_t&=x_0+\int_{0}^{t}f(s,x_s,u_s,v_s)ds+\int_{0}^{t}\int_E\gamma(s,e,x_{s-})
    \beta(s,e,x_{s-},u_s,v_s)\lambda(de)ds\\
    &\qquad\qquad\qquad \qquad\qquad \qquad\qquad+\int_{0}^{t}\sigma(s,x_s)dB_s+
            \int_{0}^{t}\int_E\gamma(s,e,x_{s-})\widetilde\mu(ds,de),\,\,t\in [0, T].
\end{array}$$
The agent $c_1$ (resp. $c_2$) controls the system with the help of
the process $u$ (resp. $v$) up to the time when she decides to stop
controlling at  $\tau$ (resp. $\sigma$), a stopping time. Then the
control of the system is stopped at $\tau\wedge \sigma$, that is to
say, when one of the agents decides first to stop controlling. As
noticed above, the advantages of the
 agents are antagonistic,
$i.e.$, there is a payoff $J(u,\tau;v,\sigma)$ between them which is
a cost (resp. a reward) for $c_1$ (resp. $c_2$). The payoff depends
on the process $(x_t)_{t\leq T}$ and is the sum of two parts, an
instantaneous and terminal payoffs (see (\ref{cout}) for its
definition). Therefore the agent $c_1$ aims at minimizing
$J(u,\tau;v,\sigma)$ while $c_2$ aims at maximizing the same payoff.
In the particular case of agents who have non control actions, the
mixed game is just the well known Dynkin game which is studied by
several authors (see e.g. \cite{laraki, LM, touzivieille} and the
references therein). Also in this paper we show that this game has a
value, $i.e.$,  the following relation holds true:
$$
\inf_{(u,\tau)}\sup_{(v,\sigma)}J(u,\tau;v,\sigma)=\sup_{(v,\sigma)}\inf_{(u,\tau)}J(u,\tau;v,\sigma).$$
The value of the game is expressed by means of a solution of a BSDE
with two reflecting barriers with a specific generator.

In the case when the filtration is Brownian ($i.e.$ the process
$(x_t)_{t\leq T}$ has no jumps), the zero-sum mixed differential
game is completely solved in \cite{[H2]} in its general setting.
However according to our knowledge the problem of zero-sum mixed
differential-integral game still open. Therefore our work completes
and closes this problem of zero-sum stochastic games of diffusion
processes with jumps.

In a financial market zero-sum games are related to recallable
options and convertible bonds. Recallable options (or israeli
options in Kifer's terminology \cite{yk}) are American options where
the issuer of the option has also the right to recall it if she
accepts to pay at least the value of the option in the market.
Therefore we have a zero-sum game between the issuer and the holder
of the option (see $e.g.$ \cite{[H2], kk, yk} for more details on
this subject). A convertible bond is a financial instrument, in
general issued by firms, with the following provisions: it pays a
fixed amount at maturity like a bond and pays coupons ; it can be
converted by the bondholder for stock or can be called by the firm.
Therefore as a game option, this makes also a zero-sum game between
the issuer and the bondholder (see $e.g.$ \cite{afv, gk, sish} and
the references therein for the literature on convertible bonds).
Also another problem that can motivate the mixed zero-sum game we
consider is the pricing of American game options or convertible
bonds under Knightian uncertainty (see e.g. \cite{knight}) with or
without defautable risk of the underlyings \cite{{bcjr},{bcjr2}}. We
will come back to this topic in a forthcoming paper.

This paper is organized as follows. In Section 2, we formulate the
problem and we recall some results related to BSDEs with one
reflecting discontinuous $rcll$ barrier. In Section 3, we introduce
the increasing and decreasing penalization schemes and we prove
their convergence. Later we show that the limits of those schemes
are the same and provides the so-called \it{local solution} for the
two barrier reflected BSDE. In Section 4 we give the main result of
this paper (Theorem \ref{thmprince}), where we establish the
existence and uniqueness of the solution of (\ref{bsdeintro}) when
the obstacles and their left limits are completely separated. We
first begin to consider the case when $f$ does not depend on
$(y,z,v)$ and in using results of Section 3 (Theorem \ref{thmimp1})
we show that existence/uniqueness, then we switch to the case where
$f$ depends only on $y$ and we use a fixed point argument to state
the existence of a solution for (\ref{bsdeintro}) (Proposition
\ref{propexist}), finally we deal with the general case. At the end,
in Section 5 we solve the mixed zero-sum differential-integral game
problem as an application of our study. $\Box$
\section{Setting of the problem and preliminary results}
Throughout this paper, $(\Omega, {\cal F},({\cal F}_t)_{t\leq T},
P)$ is a stochastic basis such that ${\cal F}_0$ contains all $\
P$-null sets of $\cal F$ and ${\cal
F}_{t+}:=\bigcap_{\epsilon>0}{\cal F}_{t+\epsilon}={\cal F}_t$,
$\forall t<T$. Moreover we assume that the filtration is generated
by the following two
mutually independent processes:\\
- a $d$-dimensional Brownian motion $(B_t)_{t\leq T},$\\
- a Poisson random measure $\mu$ on $\it R^{+}\times\it E$, where
$\it E:=R^{l}\backslash \{0\}$ $(l\geq 1$) is equipped with its
Borel $\sigma$-algebra $\mathcal {E}$, with compensator
$\nu(dt,de)=dt\lambda(de)$, such that ${\widetilde \mu([0,t]\times
A)=(\mu-\nu)([0,t]\times A)}_{t\leq T}$
 is a martingale for every $A\in \mathcal {E}$ satisfying $\lambda(A)<\infty$. The measure $\lambda$ is assumed to be a $\sigma$-finite
 on $(E,\mathcal E)$ and integrates the function $(1\wedge
 |e|^2)_{e\in E}$.
 \medskip
 Besides let us define:
\medskip

- ${\cal P}$ (resp. ${\cal P}^d$) the $\sigma$-algebra of ${\cal
F}_t$-progressively measurable (resp. predictable) sets on
$[0,T]\times \Omega$;

- ${\cal H}^{k}$ $(k\geq 1)$ the set of $\cal P$-measurable
processes $Z=(Z_t)_{t\leq T}$ with values in $R^k$ such that
$P-a.s.$, $\int_0^T|Z_s(\omega)|^2ds<\infty$ ; ${\cal H}^{2,k}$ is
the subset of the set of ${\cal H}^{k}$
 of processes $Z=(Z_t)_{t\leq T}$ $dt\otimes dP$-square integrable ;

- ${\cal S}^2$ the set of ${\cal F}_t$-adapted $rcll$ processes
$Y=(Y_t)_{t\leq T}$ such that $E[\sup_{t\leq T}|Y_t|^2]<\infty$;

-${\cal L}$ the set of ${\cal P}^d\bigotimes \mathcal E$-measurable
mappings $V:\Omega \times [0,T]\times E\rightarrow R$ such that
$P-a.s.$,
\\ $\int_{0}^{T}ds\int_E(V_s(\omega,e))^2\lambda(de)<\infty$ ; ${\cal
L}^2$ is the subset of ${\cal L}$ which contains the mappings
$V(t,\omega,e)$ which are $dt \times dP\times d\lambda$-square
integrable ;

- ${\cal A}$ the set of ${\cal P}^d$-measurable, $rcll$
non-decreasing processes $K=(K_t)_{t\leq T}$ such that $K_0=0$ and
$P-a.s.$, $K_T<\infty$ ; we denote by ${\cal A}^{2}$ the subset of
${\cal A}$ which contains processes $K$ such that $E[K_T^2]<\infty$
and by ${\cal A}^{2,c}$ the subset of ${\cal A}^2$ which contains
only continuous processes ;

- for $\pi=(\pi_t)_{t\leq T}\in {\cal S}^2$,
${\pi_-}:=({\pi}_{t-})_{t\leq T}$ is the process of its left limits,
$i.e.$, $\forall t>0, {\pi}_{t-}=lim_{s\nearrow t}\pi_s$
($\pi_{0-}=\pi_0$). On the other hand, we denote by
$\Delta\pi_t=\pi_t-\pi_{t-}$ the size of the jump of $\pi$ at $t$ ;

- a stopping time $\tau$ is called predictable if there exists a
sequence $(\tau_n)_{n\geq 0}$ of stopping times such that
$\tau_n\leq \tau$ that are strictly smaller than $\tau$ on $\{\tau
>0\}$ and increase to $\tau$ everywhere ; a stopping time $\zeta$ is
called completely inaccessible if for any predictable stopping time
$\tau$, $P[\tau=\zeta]=0$ ; the set of ${\cal F}_t$-stopping times
$\varsigma$ which take their values in $[t,T]$ is denoted by
 ${\cal T}_t$. $\Box$
\medskip

We are now given four objects: \medskip

\no $(i)$ a function ${f}:$ $(t,\omega,y,z,v)\in [0,T]\times \Omega
\times R^{1+d}\times L^2(E,\mathcal E,\lambda;\it R)\longmapsto$
$f(t,\omega,y,z,v)\in R$ such that $({f}(t,\omega,y,z,v))_{t\leq T}$
is ${\cal P}$-measurable for any $(y,z,v)\in R^{1+d}\times
L^2(E,\mathcal E,\lambda;\it R))$ and $(f(t,\omega,0,0,0))_{t\leq
T}$ belongs to ${\cal H}^{2,1}$. Moreover we assume that $f$ is
uniformly Lipschitz with respect to $(y,z,v)$, $i.e.$, there exists
a constant $C_f$ (when there is no ambiguity we omit $f$ at the
index) such that:$$ P-a.s.\,\,\, |f(t,y,z,v)-f(t,y',z',v')|\leq
C_f(|y-y'|+|z-z'|+||v-v'||), \mbox{ for any }t,y,y',z,z',v \mbox{
and }v'$$

\no $(ii)$ a random variable $\xi$ which belongs to $L^2(\Omega,
{\cal F}_T,dP)$

\no $(iii)$ two barriers $L:=(L_t)_{t\leq T}$ and $U:=(U_t)_{t\leq}$
processes of ${\cal S}^2$ which satisfy: $$ P-a.s., \forall t\leq T,
\,\, L_t\leq U_t \mbox{ and } L_T\leq \xi \leq U_T.$$

A solution for the BSDE, driven by the Brownian motion $B$ and the
independent Poisson random measure $\mu$, with two reflecting $rcll$
barriers associated with $(f,\xi,L,U)$ is a quintuple \\
$(Y,Z,V,K^+,K^-):=(Y_t,Z_t,V_t,K^+_t,K^-_t)_{t\leq T}$ of processes
with values in $R^{1+d}\times L^2_R(E,\mathcal E,\lambda)$$\times
R^{1+1}$ such that: $\forall
t\leq T$,\\
\begin{equation}\label{bsde}
      \left\{\begin{array}{l} (i)\,\, Y\in {\cal S}^2,K^\pm\in {\cal A},
      Z\in {\cal H}^{d} \mbox{ and } V\in {\cal L}\\
      (ii)\,Y_t=\xi+\int_t^{T}f(s,Y_s,Z_s,V_s)ds+(K^+_T-K^+_t)-\\\qquad\qquad\qquad\qquad\qquad\qquad
      (K^-_T-K^-_t)-\int_t^{T}Z_sdB_s -
   \int_t^{T}\int_EV_s(e)\widetilde{\mu}(ds,de),\forall t\leq T\\
   (iii)\,L\leq Y\leq U \mbox{ and if }K^{c,\pm} \mbox{ is the
   continuous part of }K^\pm \mbox{ then
   }\\\qquad\int_0^T(Y_t-L_t)dK^{c,+}_t=\int_0^T(U_t-Y_t)dK^{c,-}_t=0\\
   (iv)\,\mbox{if }K^{d,\pm} \mbox{ is the purely discontinuous part
   of }K^\pm \mbox{ then }K^{d,\pm} \mbox{ is predictable and }\\
   \qquad K^{d,-}_t=\sum\limits_{0<s\leq t}(Y_s-U_{s-})^+1_{[\Delta
                                  U_s>0]} \mbox{ and }K^{d,+}_t=\sum\limits_{0<s\leq
                                  t}(Y_s-L_{s-})^-1_{[\Delta L_s<0]}\,\,;
              \end{array}
      \right.
\end{equation}
here $x^+=max\{x,0\}$ and $x^-=-min\{x,0\}$ for any $x\in R$.
\medskip

First let us notice that obviously for arbitrary barriers $L$ and
$U$ this equation does not have a solution. Actually, if for
example, $L$ and $U$ coincide and $L$ is not a semimartingale then
we cannot find a semimartingale which equals to $L$. However as
pointed out in the introduction, under Mokobodski's condition which
reads as:
\begin{equation}\label{mokocond}\bf{[Mk]:}\,\,\left\{\begin{array}{l}
\mbox{there exist two supermartingales of } {\cal S}^2, (h_t)_{t\leq
T}\mbox{ and }(\theta_t)_{t\leq T}\mbox{ which satisfy}
\\\qquad P-a.s.,\,\,\forall t\leq T,\,\,h_t\geq 0,\,\,\theta_t\geq 0 \mbox{
and }L_t\leq h_t-\theta_t\leq U_t,
\end{array}
\right. \end{equation} there  are several works which establish
existence/uniqueness of a solution for (\ref{bsde}) (see e.g.
\cite{hh3}). So the main objective of this work is to provide
conditions on $L$ and $U$ as general as possible and easy to verify
under which equation (\ref{bsde}) has a solution. Actually in
Theorem \ref{thmprince} below we show that if the barriers $L$ and
$U$ are completely separated then the BSDE (\ref{bsde}) associated
with $(f,\xi, L,U)$ has a unique solution. This condition is
minimal. $\Box$
\medskip

To begin with we will focus on uniqueness of the solution of
(\ref{bsde}). Then we have:

\begin{proposition}\label{unique}: The RBSDE (\ref{bsde}) has at most one
solution, i.e., if $(Y,Z,V,K^+,K^-)$ and \\$(Y',Z',V',K'^+,K'^-)$
are two solutions of (\ref{bsde}), then $Y=Y', Z=Z', V=V'$ and
$K^+-K^-=K'^+-K'^-$.
\end{proposition}
$Proof:$ Since there is a lack of integrability of the processes
$(Z,V,K^+,K^-)$ and $(Z',V',K'^+,K'^-)$, we are proceeding by
localization. Actually for $k\geq 1$ let us set:
$$\tau_k:=\mbox{inf}\{t\geq 0, \integ{0}{t}(|Z_s|^2+|Z'_s|^2)ds
+\integ{0}{t}\int_E(|V_s(e)|^2+|V'_s(e)|^2)\lambda(de)ds\geq
k\}\wedge T.$$ Then the sequence $(\tau_k)_{k\geq 0}$ is
non-decreasing, of stationary type and converges to $T$ since $P$-a.s.,\\
$\int_0^T(|Z_s(\omega)|^2+|Z'_s(\omega)|^2)ds+
\int_0^T\int_E(|V_s(\omega,e)|^2+|V'_s(\omega,e)|^2)\lambda(de)ds<\infty$.
Using now It\^{o}'s formula with $(Y-Y')^2$ on
$[t\wedge\tau_k,\tau_k]$ we get:
$$
\begin{array}{ll}
(Y_{t\wedge\tau_k}-Y'_{t\wedge\tau_k})^2+\integ{t\wedge\tau_k}{\tau_k}|Z_s-Z'_s|^2ds+\sum_{t\wedge\tau_k<s\leq
\tau_k}(\Delta(Y-Y')_s)^2\\=(Y_{\tau_k}-Y'_{\tau_k})^2+2\integ{t\wedge\tau_k}{\tau_k}(Y_{s-}-Y'_{s-})(dK^+_s-dK^-_s-dK'^+_s+dK'^-_s)
\\-2\integ{t\wedge\tau_k}{\tau_k}(Y_{s}-Y'_{s})(Z_s-Z'_s)dB_s-2\integ{t\wedge\tau_k}{\tau_k}\int_E(Y_{s-}-Y'_{s-})(V_s(e)-V'_s(e))
\widetilde\mu(ds,de).
\end{array}
$$
But $(Y_{s-}-Y'_{s-})(dK^+_s-dK^-_s-dK'^+_s+dK'^-_s)\leq 0$, then
taking expectation in the two hand-sides yields:
$$E[(Y_{t\wedge\tau_k}-Y'_{t\wedge\tau_k})^2+\integ{t\wedge\tau_k}{\tau_k}|Z_s-Z'_s|^2ds
+\integ{t\wedge\tau_k}{\tau_k}\int_E|V_s(e)-V'_s(e)|^2\lambda(de)ds]\leq
E[(Y_{\tau_k}-Y'_{\tau_k})^2].$$ Using now Fatous's Lemma and
Lebesgue dominated convergence theorem w.r.t. $k$ we obtain that
$Y=Y', Z=Z', V=V'$ and $K^+ -K^-=K'^+-K'^-$. $\Box$
\medskip

Let us now recall the following result by S.Hamad\`ene and Y.Ouknine
\cite{ho} (see also \cite{essaky}) related to BSDEs with one
reflecting $rcll$ barrier.
\begin{theorem}\label{th1}\cite{ho} : The BSDE with one reflecting $rcll$
upper barrier associated with $(f,\xi,U)$ has a unique solution,
$i.e.$, there exists a unique quadruple of processes
$(Y_t,Z_t,V_t,K_t)_{t\leq T}$ such that:
\begin{equation} \label{apro}
       \left\{\begin{array}{l}
                    (i)\,\,Y\in {\cal S}^2, Z\in {\cal H}^{2,d}, V\in {\cal L}^2 \mbox{ and }K\in {\cal A}^{2}\\
                    (ii)\,\,Y_t=\xi+\int_t^Tf(s,Y_s,Z_s,V_s)ds-(K_T-K_t)-\int_t^TZ_sdB_s-\int_t^T\int_EV_s(e)
                        \tilde{\mu}(ds,de),\,\,\forall t\leq T,\\
                    (iii)\,\, \,\,Y_t\leq U_t,\,\,\forall t\leq T,\\
                    (iv)\,\,\mbox{if }K=K^{c}+K^{d}\mbox{ where }K^{c} \mbox{ (resp.
                        }K^{d}) \mbox{ is the continuous }\mbox { (resp. purely
                        discontinuous) part}\\  \mbox{ of }K \mbox{ then } K^d \mbox{
                        is predictable,} \int_{0}^{T}(U_t-Y_t)dK^{c}_t=0 \mbox{ and
                        } \Delta K_t=(Y_t-U_{t-})^+1_{[\Delta U_t>0]}, t\leq T.
              \end{array}
       \right.
\end{equation}
\mbox{Moreover, the process} Y \mbox{ can be characterized as
follows:} $\forall t\leq T$,
$$Y_t=essinf_{\tau\geq
t}E[\integ{t}{\tau}f(s,Y_s,Z_s,V_s)ds+U_\tau1_{[\tau<T]}+\xi1_{[\tau=T]}|{\cal
F}_t].$$
\end{theorem}
\begin{remark}\label{rmqimp}$(i)$ Rewriting equation $(ii)$ forwardly we
see that the predictable jumps of the process $Y$ are positive and
they are equal to the ones of $K$. The role of $K=K^c+K^d$ is to
keep the process $Y$ below $U$ and it acts with a minimal energy.
However the actions of $K^c$ and $K^d$ are complementary and not the
same. Actually $K^d$ does act only when the process $Y$ has a
predictable jump, which occurs at a predictable positive jump points
of $U$. In that case the role of $K^d$ is to make the necessary jump
to $Y$ in order to bring it below $U$. Therefore when $Y$ has a
predictable jump we compulsory have $U_-=Y_-$ and $\Delta Y_t=\Delta
K^d_t=(Y_{t}-U_{t-})^+1_{[\Delta U_t>0]\cap[Y_{t-}=U_{t-}]}$. Now
the role of $K^c$ is also to keep $Y$ below the barrier but it does
act only when $Y$ reaches $U$ either at its continuity or at its
positive jump points. This is the meaning of
$\int_0^T(U_t-Y_t)dK^{c}_t=0$.

$(ii)$ The condition of point $(iv)$ is equivalent to
$\int_0^T(U_{s-}-Y_{s-})dK_s=0$. Actually if $(iv)$ is satisfied
then $\int_0^T(U_{s-}-Y_{s-})dK_s=\int_0^T(U_{s-}-Y_{s-})dK^c_s+
\int_0^T(U_{s-}-Y_{s-})1_{[\Delta Y_s<0]}dK^d_s=0$ because
respectively the processes $Y$ and $U$ are rcll and the jumps of $K$
are predictable and occur only when $U_{t-}=Y_{t-}$. Conversely if
$\int_0^T(U_{s-}-Y_{s-})dK_s=0$ then $\int_0^T(U_{s}-Y_{s})dK^c_s=0$
and $\int_0^T(U_{s-}-Y_{s-})dK^d_s=0$. This latter implies $\Delta
K^d_t=(Y_t-U_{t-})^+1_{[U_{t-}=Y_{t-}]}=(Y_t-U_{t-})^+1_{[U_{t-}=Y_{t-}]\cap
[\Delta U_t>0]}=(Y_t-U_{t-})^+1_{[\Delta U_t>0]}$, whence the
desired result.

$(iii)$ The process $K^d$ can also be written as: $\forall t\leq T$,
$K^d_t=\sum_{0<s\leq t}(Y_{t}-U_{t-})^+1_{[\Delta U_t>0]}$. $\Box$
\end{remark}
\begin{remark}\label{rem11} In Theorem \ref{th1}, we have given the notion
of a solution of a BSDE with one upper reflecting barrier. However
one could have given the notion of a solution for a BSDE with a
lower reflecting barrier. Actually a triple $(Y,Z,V,K)$ is a
solution for the BSDE with a lower reflecting $rcll$ barrier $L$, a
coefficient $f$ and a terminal value $\xi$ iff $(-Y,-Z,-V,K)$ is a
solution  for the BSDE with a reflecting upper  rcll barrier
associated with $(-f(t,\omega,-y,-z),-\xi,-L)$. The solution $Y$ can
also be characterized as a Snell envelope of the following form,
i.e., the lowest $rcll$ supermartingale of class $[D]$ (i.e. the set
of random variables $\{Y_\tau, \tau \in {\cal T}_0\}$ is uniformly
integrable) which dominates a given process: $\forall t\leq T$,
$$Y_t=esssup_{\tau\geq t}E[\integ{t}{\tau}f(s,Y_s,Z_s,V_s)ds+
L_\tau1_{[\tau<T]}+\xi1_{[\tau=T]}|{\cal F}_t].\,\, \Box$$
\end{remark}

We will now provide a comparison result between solutions of one
barrier reflected BSDEs which plays an important role in this paper.
So assume there exists another quadruple of processes\\
$(Y',Z',V',K')$ solution for the one upper barrier reflected BSDE
associated with $(f',\xi',U)$. Then we have:
\begin{theorem}\label{comparison}
Assume that:\\
 (i) $f$ is independent of $v$\\
 (ii) $P$-a.s. for any $t\leq T$, $f(t,Y_t',Z_t')\leq
 f'(t,Y_t',Z_t',V_t')$ and $\xi\leq\xi'$.\\
Then $P$-a.s., $\forall t\leq T$, $Y_t\leq Y_t'$. Additionally, if
$f'$ does not depend on $v$ then we have also $K_t-K_s\leq
K'_t-K'_s$, for any $0\leq s\leq t\leq T$.
\end{theorem}
$Proof$: The main idea is to make use of Meyer-It\^o's formula with
$\psi(x)=(x^+)^2$, $x\in R$, and $Y-Y'$ (see e.g. \cite{protter},
pp. 221) which, after taking expectation in both hand-sides, yields:
$$
\begin{array}{ll}
{}&E[\psi(Y_t-Y'_t)+\int_{t}^{T}
1_{[Y_{s-}-Y'_{s-}>0]}|Z_s-Z'_s|^2ds \\{}&\qquad \qquad +
\somme{t<s\leq T}{}\{{\psi(Y_s-Y'_s)}-{\psi(Y_{s-}-Y'_{s-})}
-\psi'{(Y_{s-}-Y'_{s-})}\Delta (Y-Y')_s\}]\\ {}&\qquad=
E[\int_{]t,T]}{}
\psi'(Y_{s-}-Y'_{s-})\{(f(s,Y_s,Z_s)-f'(s,Y'_s,Z'_s,V'_s))ds
-d(K_s-K'_s)]\\{}&\qquad\leq E[\int_{]t,T]}
\psi'(Y_{s-}-Y'_{s-})\{(f(s,Y_s,Z_s)-f(s,Y'_s,Z'_s))ds
-d(K_s-K'_s)].
\end{array}
$$
But for any $t\leq T$, $\int_{]t,T]}
\psi'(Y_{s-}-Y'_{s-})d(K_s-K'_s)\geq 0$ since $\int_{]t,T]}
\psi'(Y_{s-}-Y'_{s-})dK'_s= \int_{]t,T]}
\psi'(Y_{s-}-Y'_{s-})\{dK^{'c}_s+dK^{'d}_s\}=0$. Actually the first
term is null since when $K^{'c}$ increases then we compulsory have
$Y'=U$ which implies that $\psi'(Y_{t-}-Y'_{t-})=0$ because $Y\leq
U$. The second term is also null because when the purely
discontinuous $K^{'d}$ increases at $t$ we should have
$Y'_{t-}=U_{t-}$ and then once more $\psi'(Y_{t-}-Y'_{t-})=0$.
Therefore for any $t\leq T$ we obtain:
$$
\begin{array}{l}
E[\psi(Y_t-Y'_t)+\int_{t}^{T} 1_{[Y_{s-}-Y'_{s-}>0]}|Z_s-Z'_s|^2ds]
\leq E[\int_t^T
\psi'(Y_{s-}-Y'_{s-})(f(s,Y_s,Z_s)-f(s,Y'_s,Z'_s))ds].
\end{array}
$$
Making use now of classical arguments to deduce that
$\psi(Y_t-Y'_t)=0$ for any $t\leq T$ and then $Y\leq Y'$.

Assume moreover now that $f'$ does not depend on $v$. In that case
the solutions of the BSDEs associated with $(f,\xi,U)$ and
$(f',\xi',U')$ respectively  can be constructed in using the
following penalization schemes. Actually for $n\geq 0$ let
$(Y^n,Z^n,V^n)$ and $(Y'^n,Z'^n,V'^n)$ defined as follows: $\forall
t\leq T$,
$$\begin{array}{l}
Y^n_t=\xi+\int_t^Tf(s,Y^n_s,Z^n_s)ds-\int_t^Tn(Y^n_s-U_s)^+ds-
\int_t^TZ^n_sdB_s-\int_t^T\int_EV^n_s(e) \tilde{\mu}(ds,de)\\\mbox{
and }\\
Y'^n_t=\xi'+\int_t^Tf'(s,Y'^n_s,Z'^n_s)ds-\int_t^Tn(Y'^n_s-U_s)^+ds
-\int_t^TZ'^n_sdB_s-\int_t^T\int_EV'^n_s(e)\tilde{\mu}(ds,de).\end{array}$$
First not that through comparison we have $Y^n\leq {Y'}^n$ for any
$n\geq 0$. On the other hand, it has been shown in (\cite{essaky},
Theorem 5.1) that the sequences $(Z^n)_{n\geq 0}$ and $(V^n)_{n\geq
0}$ (resp. $(Z'^n)_{n\geq 0}$ and $(V'^n)_{n\geq 0}$) converge in
$L^p([0,T]\times\Omega, dt\otimes dP)$ and
$L^p([0,T]\times\Omega\times U, dt\otimes dP\times d\lambda)$ to the
processes $Z$ and $V$ (resp. $Z'$ and $V'$) for any $p\in [0,2[$
(see also S.Peng \cite{[P]} in the case of Brownian filtration).
Moreover for any stopping time $\tau$ the sequence
$(Y^n_\tau)_{n\geq 1}$ and $(Y'^n_\tau)_{n\geq 1}$ converge
decreasingly to $Y_\tau$ and $Y'_\tau$ P-$a.s.$. Therefore, at least
after extracting a subsequence,  the sequences
$\int_0^\tau(Y^n_s-U_s)^+ds$ and $\int_0^\tau(Y'^n_s-U_s)^+ds$
converge in $L^p(dP)$ to $K_\tau$ and $K'_\tau$ $(p\in [0,2[)$.
Henceforth for any $s\leq t$ we have:
$$K_t-K_s=
\lim_{n\rightarrow\infty}\int_{s}^{t}n(Y_s^n-U_s)^+ds\leq
\lim_{n\rightarrow\infty}\int_{s}^{t}n({Y'}_s^n-U_s)^+ds=K'_t-K'_s
$$
since $Y^n\leq {Y'}^n$. The proof is complete.$\Box$
\begin{remark}\label{rem1} $(i)$ Using Remark \ref{rmqimp}-$(iii)$, since $Y\le Y'$ then we
obviously have  $P-a.s.$, for any $s\leq t$, $K^d_t-K^d_s\leq
K'^d_t-K'^d_s$.

\noindent (ii) If the barriers are not the same, as it is assumed in
the previous theorem, we can still get the comparison result of the
$Y's$, but the comparison of the $K's$ could fail. $\Box$
\end{remark}

Finally recall the following result related to indistinguishability
of two optional or predictable processes which is used several times
later. Let $\cal O$ be the optional $\sigma$-field on $(\Omega,{\cal
F},({\cal F}_t)_{t\leq T},P)$, $i.e.$, the $\sigma$-field generated
by the ${\cal F}_t$-adapted $rcll$ processes and $X$, $X'$ two
stochastic processes. Then we have:

\begin{theorem}\label{thmsect}(\cite{[DM]}, pp.220) Assume that for
any stopping time (resp. predictable stopping time) $\tau$ we have
P-$a.s.$, $X_\tau=X'_\tau$ and the processes $X$ and $X'$ are $\cal
O$-measurable (${\cal P}^d$-measurable). Then the processes $X$ and
$X'$ are undistinguishable. $\Box$
\end{theorem}
\section{Local solutions of BSDEs with two general $rcll$ reflecting barriers}
We are now going to show the existence of a process $Y$ which
satisfies locally the BSDE (\ref{bsde}), $i.e.$, for any stopping
time $\tau$ one can find another greater stopping time $\theta_\tau$
such that on $[\tau, \theta_\tau]$, $Y$ satisfies the BSDE
(\ref{bsde}) with terminal condition $Y_{\theta_\tau}$. The process
$Y$ will be constructed as the limit of solutions of a penalization
scheme.

For BSDEs driven by a Brownian and Poisson measure, the comparison
result between solutions does not hold in the general case,
especially when the generators depend on $v$ (see a counter-example
in \cite{[BBP]}). Therefore, we first assume that the map $f$ does
not depend on $v$, and for the sake of simplicity, we will assume
that $ f(t,\omega,y,z,v)\equiv g(t,\omega)$.

Let us now begin to analyze the increasing penalization scheme.
\subsection{The increasing penalization scheme}
Let us introduce the following increasing penalization scheme. For
$n\geq 1$, let $(Y^n_t,Z^n_t,V^n_t,K^{n}_t)_{t\leq T}$ be the
quadruple of processes with values in $R^{1+d}\times L^2(E,\cal
E,\lambda;\it R)$$\times R$ such that:
\begin{equation}
\label{sdec} \left\{
\begin{array}{l}
(i)\,\,Y^n\in{\cal S}^2, Z^n\in {\cal H}^{2,d},V^n\in {\cal L}^2
\mbox{ and }K^{n}\in {\cal A}^{2}\\
(ii) \,\,Y^n_t=\xi+\int_{t}^{T}\{g(s)+n(L_s-Y^n_s)^+\}ds
-(K_T^{n}-K_t^{n}) - \int_{t}^{T} Z^{n}_sdB_s-
\int_t^T\int_EV_s^n(e)\tilde{\mu}(ds,de)\\
(iii)\, Y^n\leq U\\
(iv)\,\,\mbox{if } K^{n,c}\mbox{ (resp. }K^{n,d}) \mbox{ is the
continuous (resp. purely discontinuous) part of }
K^{n},\,\mbox{i.e., } \\K^n=K^{n,c}+K^{n,d},\mbox{ then }
\int_{0}^{T} (U_s-Y^n_s)dK_s^{n,c}=0 \mbox{ and }K^{n,d} \mbox{ is
predictable and satisfies }\\K_t^{n,d}=\sum_{0<s\leq
t}(Y^n_s-U_{s-})^+, \forall t\leq T.
\end{array}
\right.
\end{equation}
The existence of the quadruple $(Y^n,Z^n,V^n,K^{n,-})$ is due to
Theorem \ref{th1}. Now the comparison result given in Theorem
\ref{comparison} implies that for any $n\geq 0$ we have $Y^n\leq
Y^{n+1}\leq U$ (this is the reason for which the scheme is termed as
of increasing type). Therefore there exists a right lower
semi-continuous process $Y=(Y_t)_{t\leq T}$ such that P-$a.s.$, for
any $t\leq T$, $Y_t=\lim_{n\rightarrow \infty}Y^n_t$ and $Y_t\leq
U_t$. Additionally and obviously the sequence of processes
$(Y^n)_{n\geq 0}$ converges to $Y$ in ${\cal H}^{2,1}$.

Next for an arbitrary stopping time $\tau$, let us set: $$
\begin{array}{ll}
\delta^n_\tau&:=\inf\{s\geq \tau, K^n_s-K^n_{\tau}>0\}\wedge T\\&=
\inf\{s\geq \tau, K^{n,d}_s-K^{n,d}_\tau>0\}\wedge\inf\{s\geq \tau,
K_s^{n,c}-K_\tau^{n,c}>0\}\wedge T. \end{array}$$Once more from the
comparison theorem (\ref{comparison}), $K_t^{n}-K_{\tau}^{n}\leq
K_t^{n+1}-K_{\tau}^{n+1}$, therefore $(\delta^n_\tau)_{n\geq 0}$ is
a decreasing sequence of stopping times and converges to
$\delta_\tau := \lim_{n\rightarrow\infty}\delta^n_\tau$, which is
also a stopping time. Besides note that for any $t\in
[\tau,\delta_\tau[, K^{n,d}_t-K^{n,d}_\tau=0$ for any $n\geq 0$.
\medskip

The processes $Y$ satisfies:
\begin{proposition}\label{prop1}:
For any stopping time $\tau$ it holds true:
$$P-a.s.,\,\,1_{[\delta_\tau <T]}Y_{\delta_\tau}\geq 1_{[\delta_\tau
<T]}( U_{\delta_\tau}-1_{[\delta_\tau>\tau]}(\triangle
U_{\delta_\tau})^+).$$
\end{proposition}
$Proof$: By definition of $\delta^n_\tau$,
$K^{n,c}_{\delta^n_\tau}=K^{n,c}_\tau$, hence from (\ref{sdec}), we
get that: $\forall t\in [\tau,\delta^n_\tau]$,
\begin{equation}
\label{apre}
\begin{array}{ll}
Y^n_t&= Y^n_{\delta^n_\tau}+\int_{t}^{\delta^n_\tau}\{g(s)
+n(L_s-Y^n_s)^+\}ds
-(K^{n,d}_{\delta_\tau^n}-K^{n,d}_{t})-\int_{t}^{\delta^n_\tau}
Z^{n}_s dB_s-\int_t^{\delta^n_\tau}\!\!\!\int_EV^n_s(e)\tilde{\mu}(ds,de).\\
\end{array}
\end{equation}
In this equation the term $K^{n,d}_{\delta_\tau^n}-K^{n,d}_{t}$
still remains because the process $K^{n,d}$ could have a jump at
$\delta^n_\tau$. Moreover we have: \begin{equation}
\label{eqk}\forall t\in [\tau,\delta_\tau^n],\,\,
K^{n,d}_{\delta_\tau^n}-K^{n,d}_{t}\leq 1_{[t<\delta_\tau^n]\cap
[Y^n_{\delta_\tau^n-}=U_{\delta_\tau^n-}]}(Y^n_{\delta_\tau^n}-U_{\delta_\tau^n-})^+\end{equation}
since the stoping time $\delta_\tau^n$ could be not predictable.
Next for any $n\geq 0$, we have $Y^0\leq Y^n\leq U$ then there
exists a constant $C$ such that $E[\sup_{t\leq T}|Y^n_{t}|^2]\leq
C$. Additionally since $f$ is Lipschitz then standard calculations
(see e.g. \cite{ho}) imply:
\begin{equation}
\label{eqfyz} \sup_{n\geq 0} E[
\int_{\tau}^{\delta_\tau^n}|Z_s^n|^2ds]+ \sup_{n\geq 0} E[
\int_{\tau}^{\delta_\tau^n}\!\!ds\!\!\int_E|V_s^n(e)|^2\lambda(de)]
< \infty.
\end{equation}
Then from (\ref{apre}) and (\ref{eqk}) we deduce that:
\begin{equation}
\label{eq11}
\begin{array}{l} Y_{\delta_\tau}^n1_{[\delta_\tau<T]} \geq
E[\{Y_{\delta^n_\tau}^n-1_{[\delta_\tau<\delta^n_\tau]}
(Y^n_{\delta^n_\tau}-U_{\delta^n_\tau-})^+\}1_{[\delta_\tau<T]}|{\cal
F}_{\delta_\tau}]- E[
\int_{\delta_\tau}^{\delta_\tau^n}|g(s)|ds|{\cal F}_{\delta_\tau}]
\end{array}
\end{equation}
because the random variable $1_{[\delta_\tau<T]}$ belongs to ${\cal
F}_{\delta_\tau}$.

But on the set $[\delta^n_\tau<T]$ it holds true that
$Y_{\delta^n_\tau}^n\geq U_{\delta^n_\tau}- 1_{[\delta^n_\tau>\tau]}
(\triangle U_{\delta^n_\tau})^+$. Actually thanks to Remark
\ref{rmqimp}-$(ii)$ on the set $[\delta^n_\tau>\tau]\cap
[\delta_\tau^n<T]$ we have either
$\{Y^n_{{\delta^n_\tau}-}=U_{{\delta^n_\tau}-}\mbox{ and }
Y^n_{{\delta^n_\tau}}>U_{{\delta^n_\tau}-}\}$ or
$Y^n_{\delta^n_\tau}=U_{\delta^n_\tau}$, hence
$Y_{\delta^n_\tau}^n\geq U_{\delta^n_\tau}\wedge
U_{{\delta^n_\tau}-}= U_{\delta^n_\tau}-(\triangle
U_{\delta^n_\tau})^+$. Now on $[\delta^n_\tau=\tau]\cap
[\delta_\tau^n<T]$, once more thanks to  \ref{rmqimp}-$(ii)$, there
exists a decreasing sequence of real numbers $(t^n_k)_{k\geq 0}$
converging to $\tau$ such that $Y^n_{t^n_k-}=U_{t^n_k-}$. Taking the
limit as $k\rightarrow \infty$ gives $Y^n_\tau\geq U_\tau$ since $U$
and $Y^n$ are $rcll$, whence the claim.
\medskip

Next going back to (\ref{eq11}) to obtain:
\begin{equation}
\label{111}
\begin{array}{ll} Y_{\delta_\tau}^n1_{[\delta_\tau<T]} &\geq
E[\{(U_{\delta^n_\tau}- 1_{[\delta^n_\tau>\tau]}(\triangle
U_{\delta^n_\tau})^+)1_{[\delta_\tau^n<T]}-1_{[\delta_\tau<\delta^n_\tau]}
(Y^n_{\delta^n_\tau}-U_{\delta^n_\tau-})^+\}1_{[\delta_\tau<T]}|{\cal
F}_{\delta_\tau}]\\{}&\qquad \qquad\qquad \qquad\qquad \qquad+E[\xi
1_{[\delta_\tau^n=T]\cap[\delta_\tau<T]}|{\cal F}_{\delta_\tau}]- E[
\int_{\delta_\tau}^{\delta_\tau^n}|g(s)|ds|{\cal F}_{\delta_\tau}].
\end{array}
\end{equation}
We now examine the terms of the right-hand side (hereafter $rhs$ for
short) of (\ref{111}). First note that in the space $L^1(dP)$, as
$n\rightarrow \infty$, $E[\xi
1_{[\delta_\tau^n=T]\cap[\delta_\tau<T]}|{\cal
F}_{\delta_\tau}]\rightarrow 0$ and from (\ref{eqfyz}) we deduce
also that $\integ{\delta_\tau}{\delta_\tau^n}|g(s)|ds \rightarrow 0$
since $\delta_\tau^n \rightarrow \delta_\tau$. On the other hand let
us set $A=\cap_{n\geq 0}[\delta_\tau<\delta^n_\tau]$. For $n$ large
enough we have:
$$
1_{[\delta_\tau<\delta^n_\tau]}
(Y^n_{\delta^n_\tau}-U_{\delta^n_\tau-})^+=1_A(Y^n_{\delta^n_\tau}-U_{\delta^n_\tau-})^+.$$
Therefore
$$
\limsup_{n\rightarrow \infty}1_{[\delta_\tau<\delta^n_\tau]}
(Y^n_{\delta^n_\tau}-U_{\delta^n_\tau-})^+= \limsup_{n\rightarrow
\infty}1_A(Y^n_{\delta^n_\tau}-U_{\delta^n_\tau-})^+\leq
1_A\limsup_{n\rightarrow\infty}(Y_{\delta^n_\tau}-U_{\delta^n_\tau-})^+=0.$$
Finally
$$
\begin{array}{ll}
\lim_{n\rightarrow\infty}[U_{\delta^n_\tau}-1_{[\delta^n_\tau>\tau]}(\triangle
U_{\delta^n_\tau})^+]& =U_{\delta_\tau}-1_A \lim_{n\rightarrow
\infty}1_{[\delta^n_\tau>\tau]}(\triangle U_{\delta^n_\tau})^+-
1_{A^c} \lim_{n\rightarrow \infty}1_{[\delta^n_\tau>\tau]}(\triangle
U_{\delta^n_\tau})^+\\
{}&=U_{\delta_\tau}-1_{A^c} \lim_{n\rightarrow
\infty}1_{[\delta^n_\tau>\tau]}(\triangle U_{\delta^n_\tau})^+=
U_{\delta_\tau}-1_{A^c}1_{[\delta_\tau>\tau]}(\Delta U_{\delta_\tau})^+ \\
{}&\geq U_{\delta_\tau}-1_{[\delta_\tau>\tau]}(\triangle
U_{\delta_\tau})^+
\end{array}
$$
and $1_{[\delta_\tau^n<T]\cap[\delta_\tau<T]}\rightarrow
1_{[\delta_\tau<T]}$ as $n \rightarrow \infty$. It follows that on
$[\delta_\tau<T]$ we have, at least after extracting a subsequence
and taking the limit,
$$ Y_{\delta_\tau} \geq
U_{\delta_\tau}-1_{[\delta_\tau>\tau]}(\triangle
U_{\delta_\tau})^+.$$ The proof is now complete. $\Box$\\
\begin{proposition}\label{prop43}:There exists a 4-uplet
$(Z',V',K'^{+},K'^{d,-})$ which in combination with the process $Y$
satisfies:
\begin{equation}
\label{eqlocal1} \left\{
\begin{array}{l}
(a)\, Z'\in {\cal H}^{2,d}, V'\in {\cal L}^2,\, K'^{+}\mbox{ and }K'^{d,-}\in {\cal A}^{2};\\
 (b)\,Y_t=Y_{\delta_\tau}
+\int_{t}^{\delta_\tau}g(s)ds-(K'^{d,-}_{\delta_\tau}-K'^{d,-}_{t})
+(K'^{+}_{\delta_\tau}-K'^{+}_t)\\\qquad \qquad \qquad \qquad \qquad
\qquad- \int_{t}^{\delta_\tau} Z'_s dB_s-
\int_{t}^{\delta_\tau}\int_EV'_s(e)\tilde\mu(ds,de), \,\,\forall
t\in [\tau, \delta_\tau]\\(c)\, \forall t\in [0,T], \,\,L_t\leq
Y_t\leq U_t\\ (d)\,K'^+_\tau=0 \mbox{ and if }K'^{c,+} \mbox{ (resp.
}K'^{d,+}) \mbox{ is the continuous (resp. purely discontinuous)
part of } K'^{+} \\\,\,\, \mbox{ then }K'^{d,+} \mbox{ is
predictable, }K'^{d,+}_t=\sum_{\tau<s\leq t}(L_{s-}-Y_s)^+,\,\forall
t\in [\tau,\delta_\tau]\mbox{ and }
\int_{\tau}^{\delta_\tau}(Y_s-L_s)dK'^{c,+}_s=0 \\(e)\,K'^{d,-}
\mbox{ is predictable and purely discontinuous, }K'^{d,-}_\tau=0,
\, K'^{d,-}_t=0\,\,\forall t\in [\tau,\delta_\tau[,\mbox{ and if}\\
\,\,\quad K'^{d,-}_{\delta_\tau}>0 \mbox{ then }
Y_{\delta_\tau-}=U_{\delta_\tau-} \mbox{ and
}K'^{d,-}_{\delta_\tau}=(Y_{\delta_\tau}-U_{\delta_\tau-})^+.
\end{array} \right.
\end{equation}
\end{proposition}
$Proof$: It will be divided into three steps.\\
\noindent \underline{\it{Step 1:}} Construction of the process
$K'^{d,-}$.

For $n\geq 0$ and $t\in [0,T]$ let us set
$\Delta^{n,d}_t:=K^{n,d}_{(t\vee \tau)\wedge
\delta_\tau}-K^{n,d}_{\tau}$.  The process $\Delta^{n,d}$ is purely
discontinuous and predictable.  We just focus on this latter
property. Actually for any inaccessible stopping time $\zeta$ we
have $\Delta^{n,d}_\zeta=0$ since $K^{n,d}$ is predictable. On the
other hand for any predictable stopping time $\eta$,
$\Delta^{n,d}_\eta=1_{[\tau<\eta]}K^{n,d}_\eta \in {\cal F}_{\eta
-}$. Therefore $\Delta^{n,d}$ is predictable (see $e.g.$ \cite{bs},
pp.5, Prop.4.5). Now from Remark \ref{rem1}-$(i)$ we get that for
any $n\geq 0$, $\Delta^{n,d}_t\leq \Delta^{n+1,d}_t , \forall t\leq
T$. On the other hand, $\forall t\in [\tau, \delta_\tau[$,
$\Delta^{n,d}_t=0$, and finally for any $t\in [\tau, \delta_\tau]$,
$\Delta^{n,d}_t \leq 1_{[t<\delta_\tau]\cap
[Y^n_{\delta_\tau-}=U_{\delta_\tau-}]}(Y^n_{\delta_\tau}-U_{\delta_\tau-})^+$.
It follows that $(\Delta^{n,d})_{n\geq 0}$ converges to a
non-decreasing purely discontinuous predictable $rcll$ process
$(K'^{d,-}_t)_{t\leq T}$ which satisfies $K'^{d,-}_\tau=0$ and for
any $t\in [\tau,\delta_\tau[, K'^{d,-}_t=0$. Suppose now that
$\omega$ is such that $K'^{d,-}_{\delta_\tau}(\omega)>0$ (which
implies that we compulsory have $\tau(\omega)<\delta_\tau(\omega)$).
Therefore there exists $n_0(\omega)$ such for any $n\geq n_0$ we
have $\Delta^{n,d}_{\delta_\tau}(\omega)>0$. Using Remark
\ref{rmqimp}-$(ii)$, it follows that for any $n\geq n_0$ we have
$Y^n_{\delta_\tau-}(\omega)=U_{\delta_\tau-}(\omega)$ and
$\Delta^{n,d}_{\delta_\tau}(\omega)=(Y^n_{\delta_\tau}-U_{\delta_\tau-})^+(\omega).$
Consequently we have also
$K'^{d,-}_{\delta_\tau}(\omega)=(Y_{\delta_\tau}-U_{\delta_\tau-})^+(\omega)$
and $Y_{\delta_\tau-}(\omega)=U_{\delta_\tau-}(\omega)$ since
$Y^n\leq Y\leq U$ and then the left limit of $Y(\omega)$ at
$\delta_\tau (\omega)$ exists. Thus we have established the claim
$(e)$. $\Box$
\medskip

\noindent \underline{\it{Step 2}}: $Y$ is $rcll$ on
$[\tau,\delta_\tau]$ and $Y\geq L$.

>From equation (\ref{apre}), since $\delta_\tau\leq \delta^n_\tau$
then we have: $\forall t\in
[\tau,\delta_\tau]$,\begin{equation}\label{eqtildey}\begin{array}{l}Y^n_{t}=Y^n_{\delta_\tau}+\int_{t}^{\delta_\tau}g(s)ds
+\int_{t}^{\delta_\tau}n(L_s-Y^n_s)^+ds-(K^{n,d}_{\delta_\tau}-K^{n,d}_{t})-
\int_{t}^{\delta_\tau}Z^{n}_s
dB_s-\int_{t}^{\delta_\tau}\int_EV^n_s(e)\widetilde{\mu}(ds,de).\end{array}\end{equation}
So if for $t\in [\tau,\delta_\tau]$ we set $\bar
Y^n_t=Y^n_{t}-\Delta^{n,d}_t=Y^n_{t}-(K^{n,d}_{t}-K^{n,d}_{\tau})+\int_{\tau}^{t}g(s)ds$
then $\bar Y^n$ satisfies:
$$\bar Y^n_t=\bar Y^n_{\delta_\tau}
+\int_{t}^{\delta_\tau}n(L_s-Y^n_s)^+ds-
\integ{t}{\delta_\tau}Z^{n}_s
dB_s-\int_{t}^{\delta_\tau}\!\int_EV^n_s(e)\widetilde{\mu}(ds,de).$$

Write this latter forwardly, we get that on $[\tau,\delta_\tau]$,
$\bar Y^n$ is a supermartingale for any $n$. Next it hods true that
$P-a.s.$, $\forall \,\,t\in [\tau,\delta_\tau]$, $\bar Y^n_t\leq
\bar Y^{n+1}_t$.

Actually if $\tau=\delta_\tau$ then the claim is obvious since $\bar
Y^n_t=Y^n_\tau$. Now if $t\in [\tau,\delta_\tau[\cap
[\tau<\delta_\tau]$, the claim is also obvious since for any $n\geq
0$, $\bar Y^n_t=Y^n_{t}+\int_{\tau}^{t}g(s)ds$ and we know that
$Y^n\leq Y^{n+1}$. Finally let us consider the case of
$t=\delta_\tau (\omega)$ when $\tau(\omega)<\delta_\tau(\omega)$.

First note that $\bar
Y^n_{\delta_\tau}=Y^n_{\delta_\tau}-(K^{n,d}_{\delta_\tau}-K^{n,d}_{\tau})+\int_{\tau}^{\delta_\tau}g(s)ds$.
So we are going to consider two cases.

\underline{Case 1}: If
$K^{n+1,d}_{\delta_\tau}(\omega)-K^{n+1,d}_{\tau}(\omega)=0$ then
thanks to comparison (see Remark \ref{rem1}-$(i)$) we have also
$K^{n,d}_{\delta_\tau}(\omega)-K^{n,d}_{\tau}(\omega)=0$, therefore
$\bar Y^n_{\delta_\tau}(\omega)=Y^n_{\delta_\tau}(\omega)\leq
Y^{n+1}_{\delta_\tau}(\omega)=\bar Y^{n+1}_{\delta_\tau}(\omega)$.

\underline{Case 2}: If
$K^{n+1,d}_{\delta_\tau}(\omega)-K^{n+1,d}_{\tau}(\omega)>0$ then
$\delta_\tau$ is a stopping time such that the pair
$(\omega,\delta_\tau(\omega))$ element of the graph of
$\delta_\tau$, $i.e.$ $[\![ \delta_\tau ]\!]$, does not belong to
the graph $[\![ \theta ]\!]:=\{(\omega,\theta(\omega)),\omega \in
\Omega\}$ of any inaccessible stopping time $\theta$. This is due to
the fact that the process $K^{n+1,d}$ is predictable and its jumping
times are exhausted by a countable set of disjunctive graphs of
predictable stopping times (see $e.g.$ \cite{[DM2]}, pp.128). Next
as $K^{n+1,d}_{\delta_\tau}(\omega)-K^{n+1,d}_{\tau}(\omega)=
(Y^{n+1}_{\delta_\tau}-U_{\delta_\tau-})^+1_{[Y^{n+1}_{\delta_\tau-}=U_{\delta_\tau-}]}(\omega)$
then $\bar
Y^{n+1}_{\delta_\tau}(\omega)=Y^{n+1}_{\delta_\tau-}(\omega)+\int_{\tau}^{\delta_\tau}g_s(\omega)ds=
U_{\delta_\tau-}(\omega)+\int_{\tau}^{\delta_\tau}g_s(\omega)ds$. So
if $K^{n,d}_{\delta_\tau}(\omega)-K^{n,d}_{\tau}(\omega)>0$ then it
is equal to
$(Y^{n}_{\delta_\tau}-U_{\delta_\tau-})^+1_{[Y^{n}_{\delta_\tau-}=U_{\delta_\tau-}]}(\omega)$
and $\bar
Y^{n}_{\delta_\tau}=Y^{n}_{\delta_\tau-}+\int_{\tau}^{\delta_\tau}g(s)ds=
U_{\delta_\tau-}+\int_{\tau}^{\delta_\tau}g(s)ds=\bar
Y^{n+1}_{\delta_\tau}$. Now if
$K^{n,d}_{\delta_\tau}(\omega)-K^{n,d}_{\tau}(\omega)=0$ then
$Y^n_{\delta_\tau}(\omega)=Y^n_{\delta_\tau-}(\omega)$ since
$\delta_\tau(\omega)$ cannot be equal to $\theta(\omega)$ for any
inaccessible stopping time $\theta$, therefore $Y^n(\omega)$ is
continuous at $\delta_\tau(\omega)$. It follows that $\bar
Y^n_{\delta_\tau}(\omega)=Y^n_{\delta_\tau-}(\omega)+\int_{\tau}^{\delta_\tau}g_s(\omega)ds\leq
U_{\delta_\tau-}(\omega)+\int_{\tau}^{\delta_\tau}g_s(\omega)ds=\bar
Y^{n+1}_{\delta_\tau}(\omega)$. Thus the sequence $(\bar Y^n)$ is
non-decreasing.
\medskip

Now for any $t\in [\tau, \delta_\tau]$, let us set $\bar
Y_t=\lim_{n\rightarrow \infty}\nearrow \bar Y^n_t$. As $\bar Y^n$ is
a supermartingale then $\bar Y$ is also a $rcll$ supermartingale on
$[\tau,\delta_\tau]$ (see e.g. \cite{[DM2]}, pp.86). But from the
definition of $\bar Y^n$ we obtain that $\bar Y_t= Y_t-K'^d_t
+\int_{\tau}^{t}g_sds$ and since $K'^{d,-}$ is $rcll$ then so is
$Y$. \medskip

We now focus on the second property. We know that:
$$Y^n_{\tau}=Y^n_{\delta_\tau}+\int_{\tau}^{\delta_\tau}g(s)ds
   +\int_{{\tau}}^{\delta_\tau}n(L_s-Y^n_s)^+ds-(K^{n,d}_{\delta_\tau}-K^{n,d}_{t})-
   \int_{\tau}^{\delta_\tau}Z^{n}_s
   dB_s-\int_{\tau}^{\delta_\tau}\int_EV^n_s(e)\tilde{\mu}(ds,de).$$
After taking expectation dividing by $n$ and letting $n\rightarrow
\infty$, we get
$E[\integ{\tau}{\delta_\tau}(L_s-Y_s^n)^+ds]\rightarrow 0$ since the
other terms in both hand-sides are bounded by $Cn^{-1}$. Therefore
when $\tau(\omega)<\delta_{\tau}(\omega)$,$\forall
t\in[\tau(\omega),\delta_{\tau}(\omega)[,\,Y_t(\omega)\geq
L_t(\omega)$ since $Y$ is $rcll$ on $[\tau,\delta_\tau]$. Finally
let us consider the case where $\tau(\omega)=\delta_\tau(\omega)$.
>From the previous proposition we have:
$$\begin{array}{ll}
         1_{[\tau=\delta_\tau]}Y_\tau&=1_{[\tau=\delta_\tau]\cap[\delta_\tau<T]}Y_{\delta_\tau}
             +1_{[\tau=\delta_\tau]\cap[\delta_\tau=T]}Y_T\\
         &\geq 1_{[\tau=\delta_\tau]\cap[\delta_\tau<T]}(U_{\delta_\tau}-1_{[\delta_\tau>\tau]}(\triangle
             U_{\delta_\tau})^+)+1_{[\tau=\delta_\tau]\cap[\delta_\tau=T]}\xi\\
         &\geq 1_{[\tau=\delta_\tau]\cap[\delta_\tau<T]}L_{\delta_\tau}+1_{[\tau=\delta_\tau]\cap[\delta_\tau=T]}L_T\\
         &=1_{[\tau=\delta_\tau]}L_\tau.
  \end{array}$$
It follows that for any $t\in[\tau,\delta_\tau],\,\,Y_t\geq L_t$.
Actually we cannot have $P[L_{\delta_\tau}>Y_{\delta_\tau}]>0$
because if so we obtain a contradiction in making the same reasoning
after replacing $\tau$ by $\delta_\tau$. Henceforth for any stopping
time $\tau$ we have
 $Y_\tau\geq L_\tau$ then, since $Y$ and $L$ are optional processes, from
 Theorem \ref{thmsect} we conclude that $P-a.s., \forall t\leq T$, $Y_t\geq
 L_t$.$ \Box$
 \medskip

\noindent \underline{\it{Step 3:}} $Y$ satisfies equation
(\ref{eqlocal1}).
\medskip

For $n\geq 0$, let us introduce the process $\tilde Y^n$ defined
by:$$\forall t\in[\tau,\delta_\tau],\tilde
Y^n_t=Y^n_t-\Delta^{n,d}_{t}=Y^n_t-(K^{n,d}_{t}-K^{n,d}_{\tau})$$
First note that for any $t\in [\tau,\delta_\tau[$,
$K^{n,d}_{t}-K^{n,d}_{\tau}=0$.
Therefore making the substitution in (\ref{eqtildey}) we obtain: $\forall t\in [\tau,\delta_\tau],$\\
$\begin{array}{ll}\tilde
Y^n_t&=Y^n_{\delta_\tau}-\Delta^{n,d}_{\delta_\tau}+\integ{t}{\delta_\tau}g(s)ds
+\int_{t}^{\delta_\tau}n(\tilde L^n_s-\tilde Y^n_s)^+ds-
\integ{t}{\delta_\tau}(Z^{n}_s
dB_s+\int_EV^n_s(e)\widetilde{\mu}(ds,de)),
\end{array}$\\
where $\tilde L^n_t:=L_t-\Delta^{n,d}_{t}.$ On the other hand, it
holds true that: $\forall t\in [\tau,\delta_\tau], \tilde
Y^n\geq\tilde Y^n\wedge\tilde L^n$ and
$\int_{\tau}^{\delta_\tau}(\tilde Y^n_s-\tilde Y^n_s\wedge\tilde
L^n_s)dK_s^n=0$, where $K_t^n=\int_{\tau}^{t}n(\tilde L_s-\tilde
Y_s^n)^+ds$. Henceforth thanks to Remark \ref{rem11}, we have:
$\forall t\in[\tau,\delta_\tau]$,
$$\begin{array}{ll}\tilde
Y^n_t&=\mbox{esssup }_{t\leq \sigma\leq
\delta_\tau}E[1_{[\sigma=\delta_\tau]}(Y^n_{\delta_\tau}-\Delta^{n,d}_{\delta_\tau})
+1_{[\sigma<\delta_\tau]}(\tilde L_{\sigma}^n\wedge\tilde
Y^n_{\sigma})+\int_{t}^{\sigma}g(s)ds|{\cal
F}_t]\\
{}&=\mbox{esssup }_{t\leq \sigma\leq
\delta_\tau}E[1_{[\sigma=\delta_\tau]}(Y^n_{\delta_\tau}-\Delta^{n,d}_{\delta_\tau})+1_{[\sigma<\delta_\tau]}(
L_{\sigma}\wedge Y^n_{\sigma})+\int_{t}^{\sigma}g(s)ds|{\cal
F}_t].\end{array}$$ Let us now consider the following BSDE: $\forall
t\in[0,\delta_\tau], $\begin{equation} \label{eqyprime} \left\{
\begin{array}{l} \tilde Y\in {\cal S}^2, \tilde Z\in {\cal H}^{2,d}, \tilde V
\in {\cal L}^2 \mbox{ and }\tilde K^{+}\in {\cal A}^{2};\\
\tilde Y_t= Y_{\delta_\tau}- K'^{d,-}_{\delta_\tau}
+\int_{t}^{\delta_\tau}g(s)ds +(\tilde K^{+}_{\delta_\tau}-\tilde
K^{+}_t) -\int_{t}^{\delta_{\tau}} \tilde Z_s
dB_s-\int_{t}^{\delta_\tau}\int_E \tilde V_s(e)\tilde{\mu}(ds,de),\\
\tilde Y_t\geq L_t-K'^{d,-}_{t}:=\tilde L_t,\,\mbox{ and }\tilde
K^{+}_t=\tilde K_t^{c,+}+\tilde K_t^{d,+} \mbox{ satisfies:
}\\
\int_{\tau}^{\delta_\tau}(\tilde Y_s-\tilde L_s)d\tilde
K_s^{c,+}=0,\, \tilde K^{d,+} \mbox{ is predictable and }\tilde
K^{d,+}_t=\sum_{0<s\leq t}(\tilde L_{s-}-\tilde Y_s)^+,\,\forall
t\in [0,\delta_\tau].
\end{array}
\right.
\end{equation}
The existence of the solution $(\tilde Y_t,\tilde Z_t,\tilde
V_t,\tilde K_t)_{t\leq \delta_\tau}$ is guaranteed by Theorem
\ref{th1} and Remark \ref{rem11}. Additionally we have the following
characterization for $\tilde Y$: $\forall t\in[\tau, \delta_\tau]$,
$$\tilde Y_t=\mbox{esssup }_{t\leq \sigma\leq
\delta_\tau}E[1_{[\sigma=\delta_\tau]}(Y_{\delta_\tau}-
K'^{d,-}_{\delta_\tau})+1_{[\sigma<\delta_\tau]}
L_{\sigma}+\int_{t}^{\sigma}g(s)ds|{\cal F}_t].$$ We are going now
to prove that $P-a.s.$ for any $t\in [\tau,\delta_\tau]$, $\tilde
Y^n_t\nearrow\tilde Y_t$.  Actually, $P-a.s.$, for any $t\in
[\tau,\delta_\tau]$ we have:
$$1_{[\tau\leq t<\delta_\tau]}L_{t}\wedge
Y^n_{t}+1_{[t=\delta_\tau]}(Y^n_{\delta_\tau}-\Delta^{n,d}_{\delta_\tau})\nearrow
1_{[\tau\leq t<\delta_\tau]} L_{t}+\wedge
1_{[t=\delta_\tau]}(Y_{\delta_\tau}-K'^{d,-}_{\delta_\tau}).$$ Note
that the increasing convergence of $(Y^n_{\delta_\tau}-\Delta
^{n,d}_{\delta_\tau})$ to $Y_{\delta_\tau}-K'^{d}_{\delta_\tau}$ is
obtained from \it{Step 2}. Using now Lemma \ref{l1} given in
Appendix we obtain that $\tilde Y^n\nearrow\tilde Y$, $i.e.$, for
any $t\in [\tau,\delta_\tau]$, $Y^n_t-\Delta^{n,d}_{t}\nearrow
\tilde Y_t$. Therefore for any $t\in [\tau,\delta_\tau], Y_t=\tilde
Y_t+K'^{d,-}_t$. Taking now into account the equation satisfied by
$\tilde Y$ we obtain: $\forall t\in [\tau,\delta_\tau]$,
\begin{equation}\label{eqlocy}\begin{array}{ll}
Y_t=Y_{\delta_\tau}-(K'^{d,-}_{\delta_\tau}-K'^{d,-}_{t})
+\int_{t}^{\delta_\tau}g(s)ds +(\tilde K^{c,+}_{\delta_\tau}-\tilde
K^{c,+}_t)+ (\tilde K^{d,+}_{\delta_\tau}-\tilde K^{d,+}_t)
-\int_t^{\delta_\tau}\tilde Z_s
dB_s\\
\hspace{5cm}- \int_{t}^{\delta_\tau}\int_E \tilde
V_s(e)\widetilde\mu(ds,de).
 \end{array}\end{equation}
Next let us set $K'^{c,+}_t=(\tilde K^{c,+}_{(t\vee \tau)\wedge
\delta_\tau}-\tilde K^{c,+}_{\tau})$, $t\leq T$ (and then
$K'^{c,+}_\tau=0$). Then the process $K'^{c,+}$ is non-decreasing
continuous and satisfies
$\int_{\tau}^{\delta_\tau}(Y_s-L_s)dK'^{c,+}_s=0$ since
$Y_t-L_t=\tilde Y_t-\tilde L_t$ for any $t\in [\tau,\delta_\tau]$.
Next we set $K'^{d,+}_t= (\tilde K^{d,+}_{(t\vee \tau)\wedge
\delta_\tau}-\tilde K^{d,+}_{\tau})$, $t\leq T$ (and then
$K'^{d,+}_\tau=0$). Then $K'^{d,+}$ is non-decreasing predictable
and purely discontinuous since $\tilde K^{d,+}$ is so. Finally for
$t\leq T$ let us set $Z'_t=\tilde Z_t1_{[\tau,\delta_\tau]}(t)$ and
$V'_t=\tilde V_t1_{[\tau,\delta_\tau]}(t)$. Therefore using equation
(\ref{eqlocy}) we obtain that the 5-uple
$(Y,Z',V',K'^{c,+},K'^{d,+},K'^{d,-})$ satisfies $(b)$. It remains
now to show property $(d)$.

Let $\eta$ be a predictable stopping time such that
$\eta<\delta_\tau$ and $\Delta K'^{d,+}_\eta>0$. Therefore $\Delta
K'^{d,+}_\eta=\Delta \tilde K^{d,+}_\eta=(\tilde L_{\eta-}-\tilde
Y_\eta)^+=(L_{\eta-}- Y_\eta)^+$ since $K'^{d,-}_t=0$ for any $t\in
[\tau,\delta_\tau[$. Suppose now that $\eta=\delta_\tau$ and $\Delta
K'^{d,+}_{\eta}>0$. Therefore thanks to $(\ref{eqlocy})$ we have
$0<\Delta K'^{d,+}_\eta=\Delta \tilde
K^{d,+}_\eta=Y_{\eta-}-Y_\eta+K'^{d,-}_\eta=\tilde
Y_{\eta-}-Y_\eta+K'^{d,-}_\eta=L_{\eta-}-Y_\eta+K'^{d,-}_\eta$.
Recall here that the Poisson part in (\ref{eqlocy}) have only
inaccessible jumps and $\eta$ is predictable. But if
$K'^{d,-}_\eta>0$ then $Y_{\eta-}=U_{\eta-}$ and
$K'^{d,-}_\eta=Y_{\eta}-U_{\eta-}$, then $0<\Delta
K'^{d,+}_\eta=\Delta \tilde
K^{d,+}_\eta=L_{\eta-}-Y_\eta+Y_{\eta}-U_{\eta-}\leq 0$ which is
contradictory. It follows that $K'^{d,-}_\eta=0$ and then $\Delta
K'^{d,+}_\eta=L_{\eta-}-Y_\eta=(L_{\eta-}-Y_\eta)^+.$ The proof is
now complete. $\Box$
\subsection{Analysis of the decreasing penalization scheme}
We now consider the following decreasing penalization scheme:
\begin{equation}
      \left\{\begin{array}{l}
                   (i)\,\, Y^{'n}\in {\cal S}^2,\,\,Z^{'n}\in
                       {\cal H}^{2,d},\,\,V^{'n}\in {\cal L}^2,\,\, K^{'n}\in {\cal A}^{2}\\
                   (ii)\,\,Y^{'n}_t=\xi+\integ{t}{T}\{g(s)-n(Y^{'n}_s-U_s)^+\}ds+(K^{'n}_T-K^{'n}_t)\\
                       \qquad \qquad\qquad\qquad- \integ{t}{T} Z^{'n}_sdB_s
                       -\integ{t}{T}\int_EV^{'n}_s(e)\tilde\mu(ds,de), \,\,\forall t\in[0, T]\\
                   (iii)\, Y^n\geq L\\
                                   (iv)\,\,\mbox{if } K^{'n,c}\mbox{ (resp. }K^{'n,d}) \mbox{ is the
continuous (resp. purely discontinuous) part of }
K^{'n},\,\mbox{$i.e.$, } \\K^{'n}=K^{'n,c}+K^{'n,d},\mbox{ then }
\int_{0}^{T} (Y^{'n}_s-L_{s-})dK_s^{'n,c}=0 \mbox{ and }K^{'n,d}
\mbox{ is predictable and satisfies }\\ K^{'n,d}_t=\sum_{0<s\leq
t}(L_{s-}-Y^{'n}_s)^+, \forall t\leq T.
             \end{array}
      \right.
\end{equation}
For any $n\geq 0$, the quadruple $(Y^{'n},Z^{'n},V^{'n},K^{'n})$
exists through Theorem \ref{th1}. Using once more the comparison
result Theorem \ref{comparison}, we have for any $n\geq 0$ P-$a.s.$,
$L\leq Y^{'n+1} \leq Y^{'n}$ therefore there exists a process
$Y':=(Y'_t)_{t\leq T}$ such that P-$a.s.$, $Y'\geq L$ and for any
$t\leq T$, $Y'_t=\lim_{n\rightarrow \infty}Y^{'n}_t$. Additionally
thanks to the Lebesgue dominated convergence theorem the sequence
$(Y^{'n})_{n\geq 0}$ converges to $Y'$ in ${\cal H}^{2,1}$.
\medskip
Next for any  stopping time $\tau$ and $n\geq 0$, let us set:
\begin{equation}\label{eqta2}
\begin{array}{ll} \lambda^n_\tau &:=\inf\{s\geq \tau, K^{'n}_s-K^{'n}_{\tau}>0\}\wedge T\\&= \inf{\{s\geq \tau:
K_s^{'n,d}-K_{\tau}^{'n,d}>0\}} \wedge\inf\{s\geq \tau,
K_s^{'n,c}-K_{\tau}^{'n,c}>0\}\wedge T.
\end{array}
\end{equation}
The same analysis reveals that $(\lambda^n_\tau)_{n\geq 0}$ is a
non-decreasing sequence of stopping times and converges to another
stopping time $\lambda_\tau:=\lim_{n\rightarrow \infty}
\lambda^n_\tau$. The following properties related to $Y'$, which are
the analogous of the ones of Proposition \ref{prop1} $\&$
\ref{prop43}, hold true:
\begin{proposition}:\label{prop2}${}$ $(i)$ P-a.s., $1_{[\lambda_\tau <T]}Y'_{\lambda_\tau}\leq
1_{[\lambda_\tau <T]}( L_{\lambda_\tau}+
1_{[\lambda_\tau>\tau]}(\triangle L_{\lambda_\tau})^-)$.\\
$(ii)$ There exists a 4-uplet of processes $(Z^{"}, V^{"},
K^{",-},K^{"d,+})$ which in association with $Y'$ satisfies:
\begin{equation}
\label{eq22}\! \left\{
\begin{array}{l}
\!\!\!(a)\,\,(Z^{"}, V^{"}, K^{",-},K^{"d,+})\in {\cal
H}^{2,d}\times {\mathcal
L}^2\times {\mathcal A}^{2}\times {\mathcal A}^{2}\\
\!\!\!(b) \;
Y^{'}_t=Y^{'}_{\lambda_\tau}+\int_t^{\lambda_\tau}g(s)ds
-(K^{",-}_{\lambda_\tau}-K^{",-}_t)+(K^{"d,+}_{\lambda_\tau}-K^{"d,+}_t)\\
\qquad\qquad \qquad\qquad - \int_t^{\lambda_\tau} Z^{"}_s
dB_s-\int_{t}^{\lambda_\tau}\int_EV^{"}_s(e)\tilde\mu(ds,de),
\,\,\forall t\in [\tau, \lambda_\tau]\\
\!\!\!(c) \,\,\forall t\in [0, T], L_t\leq Y'_t\leq U_t
\\
\!\!\!(d)\,K^{",-}_\tau=0 \mbox{ and if }K^{"c,-} \mbox{(resp.
}K^{"d,-}) \mbox{is the continuous (resp. purely discontinuous) part
of } K^{",-}
\\\,\,\mbox{ then }K^{"d,-} \mbox{ is predictable,
}K^{"d,-}_t=\sum_{\tau<s\leq t}(Y'_s-U_{s-})^+,\,\forall t\in
[\tau,\delta_\tau]\mbox{ and }
\int_{\tau}^{\lambda_\tau}(U_s-Y'_s)dK^{"c,-}_s=0
\\\!\!\!(e)\,K^{"d,+} \mbox{ is predictable and purely discontinuous,
}K^{"d,+}_\tau=0,
\, K^{"d,+}_t=0\,\forall t\in [\tau,\lambda_\tau[,\mbox{ and if}\\
\,\,\quad K^{"d,+}_{\lambda_\tau}>0 \mbox{ then }
Y'_{\lambda_\tau-}=L_{\lambda_\tau-} \mbox{ and
}K^{"d,+}_{\lambda_\tau}=(L_{\lambda_\tau-}-Y'_{\lambda_\tau})^+.\,\,
\Box
\end{array}
\right.
\end{equation}
\end{proposition}
${Proof}$: Actually the proof is based on the results of
Propositions \ref{prop1} $\&$ \ref{prop43}. Indeed let\\
$(\tilde{Y}^n,\tilde{Z}^n,\tilde V^n,\tilde{K}^{n,+})$ be the
solution of the BSDE defined as in (\ref{sdec}) but associated
with\\
$(-g(t),-\xi,-U,-L)$. Therefore uniqueness implies that\\
$(\tilde{Y}^n,\tilde{Z}^n,\tilde
V^n,\tilde{K}^{n,+})=(-Y^{'n},-Z^{'n},-V^{'n},K^{'n,+})$. Now the
properties $(i)$-$(ii)$ are a direct consequences of the ones proved
in Proposition \ref{prop1} $\&$ \ref{prop43}. $\Box$

\begin{remark}\label{rem3}: the process $Y'$
is $rcll$ on the interval $[\tau,\lambda_\tau]$. $\Box$
\end{remark}

\subsection{Existence of the local solution}
Recall that $Y$ (resp. $Y'$) is the limit of the increasing (resp.
decreasing) approximating scheme. Really the processes $Y$ and $Y'$
are undistinguishable as we show it now.
\begin{proposition}\label{propunique}: P-$a.s.$, for any $t\leq T$,
$Y_t=Y'_t$. Additionally $Y$ is $rcll$.
\end{proposition}

\no $Proof$: First let us point out that for any $n,m\geq 0$ and all
$t\in[0,T]$ we have $Y^n_t\leq Y^{'m}_t$. Actually to prove this
claim, we just need to apply Meyer-It\^o's formula as in Theorem
\ref{comparison} with $\psi(Y^n-Y^{'m})$ where $\psi(x)=(x^+)^{2}
(x\in R$) and to remark that:
$$\begin{array}{l}
\int_{t}^{T}\psi'(Y^n_s-Y^{'m}_s)m(Y^{'m}_s-U_s)^+ds=\int_{t}^{T}\psi'(Y^n_s-Y^{'m}_s)n(L_s-Y^n_s)^+ds=0.
\end{array}
$$ Then we argue as in Theorem \ref{comparison} to obtain that for any $t\leq T$
we have $Y_t^n\leq Y_t^{'m}$. Therefore $P-a.s., \forall t\leq T,
Y_t\leq Y'_t$.

Next let $\tau$ be a stopping time and $\mu_\tau^p$ another stopping
time defined by:
$$\mu_\tau^p:=\inf{\{ s\geq \tau : Y_s\geq U_s-p^{-1} \;\hbox{or}\;
Y'_s\leq L_s+p^{-1}\}}\wedge T$$ where $p$ is a real constant $\geq
1$. First let us notice that for all $s\in[\tau,\mu^p_\tau]\cap
[\tau< \mu^p_\tau]$ and all $n$ we have:
$$ Y^n_{s-}<U_{s-}\;\hbox{and}\; Y^{'n}_{s-}>L_{s-}.$$ Therefore
for any $s\in [\tau,\mu^p_\tau]$ we have $d(K^{n}_s+K'^{n }_s)=0$.
Now using It\^o's formula with \\$(Y^{'n}_t-Y^n_t)^2e^{2(C^2+C)t}$,
$t\in [\tau, \mu^p_\tau]$, then taking expectation in both
hand-sides yield ($C:=C_f$):
\begin{equation}
\label{eqintx} E[(Y^{'n}_\tau-Y^n_\tau)^2]\leq
e^{2(C^2+C)T}E[(Y^{'n}_{\mu^p_\tau}-Y^n_{\mu^p_\tau})^2]
\end{equation}and finally taking the limit as $n\rightarrow\infty$ to
obtain:$$ E[(Y'_\tau-Y_\tau)^2]\leq
e^{2(C^2+C)T}E[(Y'_{\mu^p_\tau}-Y_{\mu^p_\tau})^2].$$ Here note that
we are not allowed to apply It\^o formula with $Y-Y'$ because we do
not know whether $Y-Y'$ is a semimartingale on $[\tau, \mu^p_\tau]$.
Next let us show that
$E[(Y'_{\mu^p_\tau}-Y_{\mu^p_\tau})^2]\rightarrow 0$ as
$p\rightarrow \infty$. First notice that
$0\leq(Y'_{\mu^p_\tau}-Y_{\mu^p_\tau})1_{[\tau<\mu^p_\tau]}\leq
\frac{1}{p}$ since $U\geq Y'\geq Y\geq L$. Let us now focus on the
case when $\tau=\mu^p_\tau$. First we have:
\begin{equation}\label{eqega}
1_{[\tau=\mu^p_\tau]}(Y'_{\mu^p_\tau}-Y_{\mu^p_\tau})=
1_{[\tau=\mu^p_\tau]\cap [\tau<\delta_\tau\wedge
\lambda_\tau]}(Y'_{\tau}-Y_{\tau})+ 1_{[\tau=\mu^p_\tau]\cap
[\tau=\delta_\tau\wedge \lambda_\tau]}(Y'_{\tau}-Y_{\tau}).
\end{equation}
Suppose that $\omega \in [\tau=\mu^p_\tau]\cap
[\tau<\delta_\tau\wedge \lambda_\tau]$. Then there exists a sequence
of real numbers $(t_k)_{k\geq 0}$ which depends on $p$ and $\omega$
such that $t_k\searrow \tau$ as $k\rightarrow \infty$ and
$Y_{t_k}\geq U_{t_k}-\frac{1}{p}$ or $Y'_{t_k}\leq
L_{t_k}+\frac{1}{p}$. So assume we have $Y_{t_k}\geq
U_{t_k}-\frac{1}{p}$. Then taking the limit as $k\rightarrow \infty$
implies that $Y_\tau\geq U_\tau -\frac{1}{p}$ since $\omega \in
[\tau<\delta_\tau]$ and we know that $Y$ is $rcll$ on
$[\tau,\delta_\tau]$. It follows that $U_\tau\geq Y'_\tau \geq
Y_\tau \geq U_\tau-\frac{1}{p}$. In the same way we can show that if
$Y'_{t_k}\leq L_{t_k}+\frac{1}{p}$ then $L_\tau \leq Y_\tau\leq
Y'_\tau \leq L_\tau +\frac{1}{p}.$ Therefore
$1_{[\tau=\mu^p_\tau]\cap [\tau<\delta_\tau\wedge
\lambda_\tau]}(Y'_{\tau}-Y_{\tau})\leq \frac{1}{p}$. Finally let us
deal with the second term of (\ref{eqega}). We have:
$$\begin{array}{ll}
1_{[\tau=\delta_\tau\wedge \lambda_\tau]}(Y'_{\tau}-Y_{\tau}) &=
1_{[\tau=\delta_\tau\wedge \lambda_\tau]\cap
[\tau<T]}(Y'_{\tau}-Y_{\tau})\\&=1_{[\tau=\delta_\tau]\cap[\tau<T]\cap[\delta_\tau\leq\lambda_\tau]}
(Y'_{\delta_\tau}-Y_{\delta_\tau})
+1_{[\tau=\lambda_\tau]\cap[\tau<T]\cap[\lambda_\tau<\delta_\tau]}(Y'_{\lambda_\tau}-Y_{\lambda_\tau})
\\&=1_{[\tau=\delta_\tau]\cap[\tau<T]\cap[\delta_\tau\leq\lambda_\tau]}(Y'_{\delta_\tau}-U_{\delta_\tau})
+1_{[\tau=\lambda_\tau]\cap[\tau<T]\cap[\lambda_\tau<\delta_\tau]}(L_{\lambda_\tau}-Y_{\lambda_\tau})
\\&\leq 0
\end{array}$$ because in that case, taking into account of \ref{prop1} $\&$
\ref{prop2}-$(i)$, we have either $Y_{\delta_\tau}=U_{\delta_\tau}$
or $Y'_{\lambda_\tau}=L_{\lambda_\tau}$ and we know that $U\geq
Y'\geq Y\geq L$.

It follows that $0\leq
(Y'_{\mu^p_\tau}-Y_{\mu^p_\tau})^2=1_{[\tau<\mu^{p}_{\tau}]}(Y'_{\mu^p_\tau}-Y_{\mu^p_\tau})^2
+1_{[\tau=\mu^{p}_{\tau}]}(Y'_{\mu^p_\tau}-Y_{\mu^p_\tau})^2\leq
\frac{1}{p^2}$, then taking the limit as $p\rightarrow \infty$ in
(\ref{eqintx}) we deduce that $Y_\tau=Y'_\tau$. As $\tau$ is an
arbitrary stopping time then $P-a.s.$, $Y=Y'$.

We are now going to deal with the second property. For any $t\leq
T$, we have: $U_t\geq Y_t\geq Y_t^n$ and $L_t\leq Y_t^{'}\leq
Y_t^{'n}$, hence from the right continuity of $Y^n$ and $Y^{'n}$ we
have:
$$\liminf\limits_{s\downarrow t}Y_s\geq\liminf\limits_{s\downarrow
t}Y_s^n=Y_t^n \mbox{ and } \limsup\limits_{s\downarrow
t}Y_s=\limsup\limits_{s\downarrow t}Y_s^{'}\leq
\limsup\limits_{s\downarrow t}Y_s^{'n}=Y_t^{'n}.$$ Letting
$n\rightarrow\infty$ we get the right continuity of $Y$ since
$Y=Y'$. Let us now show that $Y$ has left limits. Define the
predictable processes $\bar Y$ and $\tilde Y$ as following: $\bar
Y_t=\liminf\limits_{s\uparrow t}Y_s$ and $\tilde
Y_t=\limsup\limits_{s\uparrow t}Y_s$. Then, we only need to prove
that for any predictable stopping time $\tau$, we have $\bar
Y_\tau=\tilde Y_\tau$. Let $(s_k)_k$ be a sequence of stopping times
which announce $\tau$. Then we have:
$$\tilde Y_\tau=\limsup\limits_{s_k\uparrow \tau}Y_{s_k}=
\limsup\limits_{s_k\uparrow
\tau}Y_{s_k}^{'}\leq\limsup\limits_{s_k\uparrow \tau}Y_{s_k}^{'n}
=\lim\limits_{s_n\uparrow
\tau}Y_{s_k}^{'n}=Y^{'n}_{\tau-}=Y^{'n}_\tau+(L_{\tau-}-Y_{\tau}^{'n})^+.$$Letting
now $n\rightarrow\infty$, we obtain, $\tilde Y_\tau\leq
Y_\tau+(L_{\tau-}-Y_\tau)^+$. Similarly, we can also get that\\
$\bar Y_\tau\geq Y_\tau-(Y_\tau-U_{\tau-})^+$. Since we obviously
have $L_{\tau-}\leq\bar Y_\tau\leq\tilde Y_\tau\leq U_{\tau-}$ then
combining the three inequalities yields: $$L_{\tau-}\vee
(Y_\tau-(Y_\tau-U_{\tau-})^+)\leq\bar Y_\tau\leq\tilde Y_\tau\leq
U_{\tau-}\wedge (Y_\tau+(L_{\tau-}-Y_\tau)^+)$$ Note that the
right-hand and the left-hand sides are equal to
$L_{\tau-}1_{[Y_\tau<L_{\tau-}]}+Y_\tau1_{[L_{\tau-}\leq Y_\tau\leq
U_{\tau-}]}\\+U_{\tau-}1_{[Y_\tau>U_{\tau-}]}$. Therefore for any
predictable stopping time $\tau$, $\tilde Y_\tau=\bar Y_\tau$, hence
due to the predictable section theorem (Theorem \ref{thmsect}),
$\tilde Y$ and $\bar Y$ are undistinguishable. It follow that
$\lim\limits_{s\nearrow t}Y_s$ exists for any $t\leq T$ and then $Y$
has left limits. $\Box$
\medskip

Through Propositions \ref{prop1}, \ref{prop43} and the previous one
we have:
\begin{corollary} \label{cor1}The process $Y$ satisfies:
$$Y_{\delta_\tau}\geq U_{\delta_\tau}-1_{[\tau<\delta_\tau]}(\Delta U_{\delta_\tau})^+ \mbox{ on } [\delta_\tau<T]
\mbox{ and } Y_{\lambda_\tau}\leq
L_{\lambda_\tau}+1_{[\tau<\lambda_\tau]}(\Delta L_{\lambda_\tau})^-
\mbox{ on } [\lambda_\tau<T]. \Box$$
\end{corollary}

Summing up now the results obtained in Propositions \ref{prop1},
\ref{prop43} and \ref{prop2}, we have the following result related
to the existence of local solutions for the BSDE (\ref{bsde}).
\begin{theorem}\label{thm3}: There exists a process
$Y:=(Y_t)_{t\in [0,T]}$ such that:

$(1)$ $Y$ is $\cal P$-measurable, $rcll$ and satisfies : $Y_T=\xi$

$(2)$ for any stopping time $\tau$ there exists a stopping time
$\theta_\tau \geq \tau$, P-$a.s.$, and a quadruple of processes
$(Z^{\tau}, V^{\tau}, K^{\tau,+}, K^{\tau, -})\in {\cal
H}^{2,d}\times {\cal L}^2\times {\mathcal A}^{2}\times {\mathcal
A}^{2}$ ($K^{\tau,\pm}_\tau=0$) such that:

\no The process $Y$ satisfies the following equation which we notice
hereafter $\cal BL$$(\xi,g,L,U)$: P-a.s.,
$$\left\{\begin{array}{l}
(i) \; Y_t=Y_{\theta_\tau}+\integ{t}{\theta_\tau}g(s)ds
+(K_{\theta_\tau}^{\tau,+}-K_{t}^{\tau,+})-(K_{\theta_\tau}^{\tau,-}-K_{t}^{\tau,-})
- \integ{t}{\theta_\tau} Z^\tau_s dB_s-\integ{t}{\theta_\tau}\int_E V^\tau_s\tilde\mu(ds,de),\\
\qquad\qquad\forall t\in [\tau, \theta_\tau ]\\ (ii)\;P-$a.s.$, \forall t\in[0,T] ,\, L_t \leq Y_t\leq U_t\\
(iii) \; \integ{\tau}{\theta_\tau} (U_s-Y_s)dK_s^{\tau
c,-}=\integ{\tau}{\theta_\tau} (Y_s-L_s)dK_s^{\tau c,+}=0, \mbox{
where }K^{\tau c,\pm} \mbox{ is the continuous part of }\\\qquad
K^{\tau,\pm} \\
(iv)\, \mbox{ the process }K^{\tau,+}\mbox{ and } K^{\tau, -} \mbox{
are predictable and } \forall t\in [\tau,\theta_\tau],\,\,K^{\tau
d,+}_t=\sum_{\tau <s\leq t}(L_{s-}-Y_s)^+ \\\qquad \mbox{ and
}K^{\tau d,-}_t=\sum_{\tau <s\leq t}(Y_s-U_{s-})^+ ,\mbox{ where }
K^{\tau d,\pm} \mbox{ is the purely discontinuous part of
}K^{\tau,\pm}.
\end{array}
\right.
$$
Hereafter we say that $Y$ is the solution of ${\cal
BL}{(g,\xi,L,U)}$.
\end{theorem}
\no $Proof$: Let $Y:=(Y_t)_{t\leq T}$ be the adapted process defined
as the limit of the increasing (or decreasing) scheme. Obviously it
is $rcll$ and satisfies, $L\leq Y\leq U$ and $Y_T=\xi$, P-$a.s.$.
\medskip

Let us now focus on $(2)$. Let $\tau$ be a stopping time, let
$\delta_\tau$ be the stopping time defined in the previous section
and finally let us set $\theta_\tau=\lambda_{\delta_\tau}$. Thanks
to Proposition \ref{prop2}, there exists \\$(Z^{"\delta_\tau},
V^{"\delta_\tau}, K^{"\delta_\tau d,+},K^{"\delta_\tau,-})$ (which
we only denote $(Z^{"},V^{"}, K^{"d,+},K^{",-})$) such that:
\begin{equation}
\label{eq22x} \!\!\!\left\{
\begin{array}{l}
\!\!\!(a)\,\,(Z^{"}, V^{"}, K^{",-},K^{"d,+})\in {\cal
H}^{2,d}\times {\mathcal
L}^2\times {\mathcal A}^{2}\times {\mathcal A}^{2}\\
\!\!\!(b) \; Y_t=Y_{\theta_\tau}+\int_t^{\theta_\tau}g(s)ds
-(K^{",-}_{\theta_\tau}-K^{",-}_t)+(K^{"d,+}_{\theta_\tau}-K^{"d,+}_t)\\
\qquad\qquad \qquad\qquad \qquad\qquad\qquad\qquad -
\int_t^{\theta_\tau} Z^{"}_s
dB_s-\int_{t}^{\theta_\tau}\int_EV^{"}_s(e)\tilde\mu(ds,de),
\,\,\forall t\in [\delta_\tau, \theta_\tau]\\
\!\!\!(c)\,K^{",-}_{\delta_\tau}=0 \mbox{ and if }K^{"c,-} \mbox{
(resp. }K^{"d,-}) \mbox{ is the continuous (resp. purely
discontinuous) part of } K^{",-}
\\\, \mbox{then }K^{"d,-} \mbox{ is predictable,
}K^{"d,-}_t=\sum_{\delta_\tau<s\leq t}(Y_s-U_{s-})^+,\,\forall t\in
[\delta_\tau,\theta_\tau]\mbox{ and }
\int_{\delta_\tau}^{\theta_\tau}(U_s-Y_s)dK^{"c,-}_s=0
\\\!\!\!(d)\,K^{"d,+} \mbox{ is predictable and purely discontinuous,
}K^{"d,+}_{\delta_\tau}=0,
\,\,\,\, K^{"d,+}_t=0 \,\forall t\in [\delta_\tau,\theta_\tau[,\mbox{ and if}\\
\,\,\quad K^{"d,+}_{\theta_\tau}>0 \mbox{ then }
Y_{\theta_\tau-}=L_{\theta_\tau-} \mbox{ and
}K^{"d,+}_{\theta_\tau}=(L_{\theta_\tau-}-Y_{\theta_\tau})^+.
\end{array}
\right.
\end{equation}
Now for any $t\leq T$, let us set:
\medskip

$(i)$ $Z^\tau_t:= Z^{'}_t 1_{[\tau\leq t\leq \delta_\tau]} +
Z_t^{''} 1_{[\delta_\tau <t\leq \theta_\tau]}$ and  $V^\tau_t:=
V^{'}_t 1_{[\tau\leq t\leq \delta_\tau]} + V_t^{''} 1_{[\delta_\tau
<t\leq \theta_\tau]}$

$(ii)$ $K^{\tau c,+}_t:=K^{'c,+}_{(t\wedge \delta_\tau)\vee \tau}$,
$K^{\tau c,-}_t:=K^{"c,-}_{(t\wedge \theta_\tau)\vee \delta_\tau}$,
$K^{\tau d,+}_t:=K^{'d,+}_{(t\wedge \delta_\tau)\vee
\tau}+K^{"d,+}_{(t\wedge \theta_\tau)\vee \delta_\tau}$, $K^{\tau
d,-}_t:=K^{'d,-}_{(t\wedge \delta_\tau)\vee \tau}+K^{"d,-}_{(t\wedge
\theta_\tau)\vee \delta_\tau}$ and finally $K^{\tau,+}=K^{\tau
c,+}+K^{\tau d,+}$ and $K^{\tau,-}=K^{\tau c,-}+K^{\tau d,-}$.
\medskip

\noindent The constructions of $Z^\tau$ and $V^\tau$ are the
concatenations of $Z^{'}$ and $Z^{''}$ (resp. $V^{'}$ and $V^{''}$).
The same happens for the construction of the processes $K^{\tau
c,\pm}$ and $K^{\tau d,\pm}$.

The process $Z^\tau$ (resp. $V^\tau$) belongs to ${\cal H}^{2,d}$
(resp. ${\cal L}^2$) and, through their definitions, the processes
$K^{\tau d,\pm}$ are non-decreasing, purely discontinuous and
predictable, $K^{\tau c,\pm}$ are non-decreasing, predictable and
continuous, finally all of them belong to ${\cal A}^{2}$.

Next let us show that $Y$, $Z^\tau$, $V^\tau$ and $K^{\tau,\pm}$
enjoy the relations of $(2)$.
\medskip

Let $t\in [\tau, \theta_\tau]$. First assume that $t\in
[\delta_\tau, \theta_\tau]$. Then from (\ref{eq22x}) and the above
definitions we have:
\begin{equation}
\label{eq23}
\begin{array}{l}
Y_{\theta_\tau}+\int_{t}^{\theta_\tau}g(s)ds
+\int_{t}^{\theta_\tau}d(K_s^{\tau,+}-K_s^{\tau,-})-
\int_{t}^{\theta_\tau} Z^\tau_s dB_s-\int_{t}^{\theta_\tau}\int_E V^\tau_s\tilde\mu(ds,de)\\
\quad =Y_{\theta_\tau}+\int_{t}^{\theta_\tau}g(s)ds
-(K^{",-}_{\theta_\tau}-K^{",-}_t)+(K^{"d,+}_{\theta_\tau}-K^{"d,+}_t)\\
\qquad\qquad-\int_{t}^{\theta_\tau}Z_s^{"}
dB_s-\int_{t}^{\theta_\tau}\int_E V^{"}_s\tilde\mu(ds,de)\\
\quad =Y_t.
\end{array}
\end{equation}
Suppose now that $t\in [\tau, \delta_\tau[$, then we have:
$$\begin{array}{ll}
Y_{\theta_\tau}&+\int_{t}^{\theta_\tau}g(s)ds
+\int_{t}^{\theta_\tau}d(K_s^{\tau,+}-K_s^{\tau,-})-
\int_{t}^{\theta_\tau} Z^\tau_s dB_s-\int_{t}^{\theta_\tau}\int_E V^\tau_s\tilde\mu(ds,de)\\
{}&=Y_{\theta_\tau}+\int_{\delta_\tau}^{\theta_\tau}g(s)ds
+\int_{\delta_\tau}^{\theta_\tau}d(K_s^{\tau,+}-K_s^{\tau,-})-
\int_{\delta_\tau}^{\theta_\tau} Z^\tau_s
dB_s-\int_{\delta_\tau}^{\theta_\tau}\int_E
V^\tau_s\tilde\mu(ds,de)\\{}& \qquad + \int_{t}^{\delta_\tau}g(s)ds
+\int_{t}^{\delta_\tau}d(K_s^{\tau,+}-K_s^{\tau,-})-
\int_{t}^{\delta_\tau} Z^\tau_s dB_s-\int_{t}^{\delta_\tau}\int_E V^\tau_s\tilde\mu(ds,de)\\
{}&=Y_{\delta_\tau}+\int_{t}^{\delta_\tau}g(s)ds
+\int_{t}^{\delta_\tau}d(K_s^{\tau,+}-K_s^{\tau,-})-
\int_{t}^{\delta_\tau} Z^\tau_s dB_s-\int_{t}^{\delta_\tau}\int_E
V^\tau_s\tilde\mu(ds,de)
\\{}&= Y_{\delta_\tau}+\int_{t}^{\delta_\tau}g(s)ds
+(K^{',+}_{\delta_\tau}-K^{',+}_t)-(K^{'d,-}_{\delta_\tau}-K^{'d,-}_t)-
\int_{t}^{\delta_\tau} Z^{'}_s dB_s-\int_{t}^{\delta_\tau}\int_E
V^{'}_s\tilde\mu(ds,de)\\{}&=Y_t.
\end{array}
$$
Therefore the processes $(Y,Z^\tau,V^\tau,K^{\tau, +},K^{\tau, -})$
satisfy equation $(2.i)$.

Next from the definitions of $K^{\tau, +}$ and $K^{\tau, -}$,
(\ref{eq22x})-$(c)$ and  (\ref{eqlocal1})-$(d)$ we have:
$$\begin{array}{l}\int_{\tau}^{\theta_\tau}(Y_s-L_s)dK^{\tau
c,+}_s=\int_{\tau}^{\delta_\tau}(Y_s-L_s)dK^{'c,+}_s=0\mbox{ and }
\int_{\tau}^{\theta_\tau}(U_s-Y_s)dK^{\tau
c,-}_s=\int_{\delta_\tau}^{\theta_\tau}(U_s-Y_s)dK^{"\tau
c,-}_s=0.\end{array}$$ Now let $\eta$ be a predictable stopping time
such that $\tau\leq \eta \leq \theta_\tau$. Therefore thanks to
relation $(2.i)$ we have:
$$
\Delta Y_\tau=\Delta K^{\tau d,-}_\eta - \Delta K^{\tau d,+}_\eta.$$
But $\{\Delta K^{\tau d,-}>0\}\subset \{Y\geq U_{-}\}$ and $\{\Delta
K^{\tau d,+}>0\}\subset \{Y\leq L_{-}\}$. As $L_{-}\leq U_{-}$ then
$\Delta K^{\tau d,-}$ and $\Delta K^{\tau d,+}_\eta$ cannot jump in
the same time. Henceforth the positive (resp. negative) predictable
jumps of $Y$ are the same as the ones of $K^{\tau d,-}$ (resp.
$K^{\tau d,+}$).

Assume now that $\Delta K^{\tau d,+}_\eta >0$. Therefore the
definitions of $K^{\tau d,+}$, $K^{'d,+}$ and $K^{"d,+}$ imply that:
$$\begin{array}{ll}
\Delta K^{\tau d,+}_\eta&= \Delta K^{'d,+}_{\eta}1_{[\tau<\eta\leq
\delta_\tau]}+\Delta K^{"d,+}_{\eta}1_{[\eta=\theta_\tau]}\\
{}&=(L_{\eta-}-Y_{\eta})^+ 1_{[\tau<\eta\leq
\delta_\tau]}+1_{[\eta=\theta_\tau]}
(L_{\theta_\tau-}-Y_{\theta_\tau})^+=(L_{\eta-}-Y_{\eta})^+
\end{array}
$$
because from (\ref{eq22x}) we deduce that on the interval
$]\delta_\tau,\theta_\tau[$ the process $Y$ does no have any
predictable negative jump. Similarly for any predictable stopping
time $\eta$ such that $\tau\leq \eta \leq \theta_\tau$ and $\Delta
K^{\tau d,-}_\eta >0$, $\Delta K^{\tau d,-}_\eta
=(Y_\eta-U_{\eta-}).$ Thus we have proved $(2.iv)$. $\Box$
\medskip
\begin{remark} \label{rem32}
When the process $Y$ is fixed, from Proposition \ref{unique} we
deduce that the quadruple \\
$(Z^{\tau},V^{\tau},K^{\tau,+},K^{\tau,-})$ is unique on
$[\tau,\theta_\tau]$. $\Box$
\end{remark}

We are now ready to show that BSDE (\ref{bsde}) has a solution. We
first focus on the case when the generator $f$ does not depend on
$(y,z,v)$ and later we deal with the general case.

\section{Existence of a global solution for the BSDE with two
completely separated $rcll$ barriers} Let us assume that the
barriers $L$ and $U$ and their left limits are completely separated,
$i.e.$, they satisfy the following assumption: \vskip 3mm
\begin{center}
\bf {[H]}: $P-a.s., \forall t\leq T$, $L_t<U_t$ and $L_{t-}<U_{t-}$.
\end{center}
Then we have:
\begin{theorem}\label{thmimp1}: Under Assumption \bf {[H]}, the BSDE
associated with $(g(t),\xi,L,U)$ has a unique solution.
\end{theorem}
\no $Proof:$ Let $Y$ be the $rcll$ process defined in Theorem
\ref{thm3}. Then for any $n\geq 1$, there exists a stopping time
$\gamma_n$, defined recursively as
$\gamma_0=0,\gamma_n=\theta_{\gamma_{n-1}}$, and a unique quadruple
$(Z^n,V^n,K^{n,+},K^{n,-})$ which belongs to ${\H}^{2,d} \times
{\L}^2\times {\mathcal A}^{2}\times {\mathcal A}^{2}$ and which with
the process $Y$ satisfy ${\cal BL}{(\xi,g,L,U)}$ on $[\gamma_{n-1},
\gamma_{n}]$.

First let us show that for any $n\geq 1$,
$P[(\gamma_{n-1}=\gamma_n)\cap (\gamma_n<T)]=0$.

Actually let $\omega$ be such that
$\gamma_{n-1}(\omega)=\gamma_n(\omega)$ and $\gamma_n(\omega)<T$.
Then using the properties of Corollary \ref{cor1}, we have
$Y_{\gamma_n}(\omega)=L_{\gamma_n}(\omega)=U_{\gamma_n}(\omega)$. As
we know that P-$a.s.$, $L<U$ then $P[(\gamma_{n-1}=\gamma_n)\cap
(\gamma_n<T)]=0$.

We will now prove that the  sequence $(\gamma_n)_{n\geq 1}$ is of
stationary type, $i.e.$, $P[\omega, \gamma_n(\omega)<T, \forall
n\geq 1]=0$. In other words for $\omega$ fixed there exists an
integer rank $n_0(\omega)$ such that for $n\geq n_0(\omega)$
$\gamma_n(\omega)=\gamma_{n+1}(\omega)=T$. Indeed let us set
$A=\cap_{n\geq 1}(\gamma_n<T)$ and let us show that $P(A)=0$. Let
$\omega\in A$ and let us set $\gamma(\omega):=\lim_{n\rightarrow
\infty}\gamma_n(\omega)$. Using once more the inequalities of
Corollary \ref{cor1}, there exist two sequences
$(t_n(\omega))_{n\geq 1}$ and $(t'_n(\omega))_{n\geq 1}$ such that
for any $n\geq 1$, $t_n, t'_n\in [\gamma_{n-1},\gamma_{n}]$,
$Y_{t_n}\geq U_{t_n}\wedge U_{t_n-}=U_{t_n}-(\Delta U_{t_n})^+$ and
$Y_{t'_n}\leq L_{t'_n}\vee L_{t'_n-}=L_{t'_n}+(\Delta L_{t'_n})^-$.
Now as $(t_n)_{n\geq 1}$ and $(t'_n)_{n\geq 1}$ are not of
stationary type since $\gamma_n(\omega)<\gamma_{n+1}(\omega)$ then
taking the limit as $n\rightarrow \infty$ to obtain that
$Y_{\gamma-}(\omega)\leq L_{\gamma-}(\omega)\leq
U_{\gamma-}(\omega)\leq Y_{\gamma-}(\omega)$. It means that the
previous inequalities are equalities and then
$L_{\gamma-}(\omega)=U_{\gamma-}(\omega)$. But this is impossible
since P-$a.s.$, $\forall \,\,t\leq T$, $L_{t-}<U_{t-}$. It follows
that $(\gamma_n)_{n\geq 1}$ is of stationary type.

Next let us introduce the following processes $Z,V,K^{\pm}$:
$P-a.s.$, for any $t\leq T$, one sets:
$$\begin{array}{ll}
Z_t=Z_t^11_{[0,\gamma_1]}(t)+\sum_{n\geq
1}Z_t^{n+1}1_{]\gamma_n,\gamma_{n+1}]},
V_t=V_t^11_{[0,\gamma_1]}(t)+\sum_{n\geq 1}V_t^{n+1}1_{]\gamma_n,\gamma_{n+1}]}\\
K_t^{c,\pm}=\left\{\begin{array}{ll}
K_t^{1c,\pm} & \mbox{if $t \in [0,\gamma_1]$}\\
K_{\gamma_{n}}^{c,\pm}+K_t^{(n+1)c,\pm} & \mbox{if $t \in
]\gamma_{n},\gamma_{n+1}]$}
                  \end{array}
          \right.\\
          K_t^{d,\pm}=\left\{\begin{array}{ll}
K_t^{1d,\pm} & \mbox{if $t \in [0,\gamma_1]$}\\
K_{\gamma_{n}}^{d,\pm}+K_t^{(n+1)d,\pm} & \mbox{if $t \in
]\gamma_{n},\gamma_{n+1}]$}.
                  \end{array}\right.\\

\end{array}$$

\noindent Then a concatenation procedure and the same analysis as
the one in Theorem 5.1 in \cite{hho} imply that the 5-uplet
$(Y,Z,V,K^{\pm})$ verify the BSDE and the uniqueness of the solution
has been shown in Proposition 2.1.

\begin{remark}\label{inte}: The sequence of stopping times $(\gamma_k)_{k\geq 0}$ will be called associated with
the solution $(Y,Z,V,K^\pm)$. Also note that for any $k$, we have
the following local integrability of the processes $Z,V$ and
$K^\pm$:
$$E[\integ{0}{\gamma_{k}}ds\{|Z_s|^2+\int_E|V_s(e)|^2\lambda(de)\}+
(K^+_{\gamma_k})^2+(K^-_{\gamma_k})^2]<\infty. \,\Box$$
\end{remark}

We are now going to investigate under which conditions Mokobodski's
condition
 introduced in (\ref{mokocond}) is verified. Actually we will show that it is locally satisfied  when $\bf {[H]}$ is fulfilled.

\begin{proposition}\label{weakm}
Under \bf{[H]},  there exists a sequence $(\gamma_k)_{k\geq 0}$ of
stopping times such that:\\
(i) for any $k\geq 0, \gamma_k\leq \gamma_{k+1}\, \mbox {and the
sequence is of stationary type, i.e.}\, P[\gamma_k<T,\forall k\geq 0]=0(\gamma_0=0)$;\\
(ii) for any $k\geq 0$, there exists a pair $(h^k,h^{'k})$ of
non-negative supermartingales which belong to ${\cal S}^2$ such
that:$$P-a.s., \forall t\leq \gamma_k, L_t\leq h^k_t-h^{'k}_t\leq
U_t.$$
\end{proposition}
$Proof:$ Let $(Y,Z,V,K^+,K^-)$ be the solution of the RBSDE
associated with $(0,\xi,L,U)$ which exists thanks to Theorem
\ref{thmimp1}. Let $(\gamma_k)_{k\geq 0}$ be the sequence of
stopping times associated with this solution (see Remark
\ref{inte}). By construction this sequence satisfies the claim
$(i)$. Let us focus on $(ii)$. For $k\geq 1$ and $t\leq T$ one sets:
$$\begin{array}{ll}h_{t\wedge \gamma_k}^k=E[Y^+_{\gamma_k}+(K^+_{\gamma_k}-K^+_{t\wedge \gamma_k})|{\cal
F}_{t\wedge \gamma_k}]\mbox{ and } h_{t\wedge
\gamma_k}^{'k}=E[Y^-_{\gamma_k}+(K^-_{\gamma_k}-K^-_{t\wedge
\gamma_k})|{\cal F}_{t\wedge \gamma_k}]\end{array}$$where
$Y^+_{\gamma_k}=\max\{Y_{\gamma_k},0\}$ and
$Y^-_{\gamma_k}=\max\{-Y_{\gamma_k},0\}$. Then $h^k,h^{'k}$ are
supermartingales of  ${\cal S}^2$ which satisfy $L_t\leq
h_t^k-h_t^{'k}\leq U_t$ for any $t\leq \gamma_k$ since
$E[\int_{0}^{\gamma_{k}}ds\{|Z_s|^2+\int_E|V_s(e)|^2\lambda(de)\}+
(K^+_{\gamma_k})^2+(K^-_{\gamma_k})^2]<\infty$. Thus we have the
desired result. $\Box$
\bigskip

Next with the help of this result we will be able to prove that the
BSDE (\ref{bsde}) has a solution in the case when the function $f$
depends also on $y$, $i.e.$, $f(t,\omega,y,z,v)=f(t,\omega,y)$.
Actually we have:

\begin{proposition}\label{propexist}
Under \bf{[H]}, the BSDE associated with $(f(t,y),\xi,L,U)$ has a
unique solution.
\end{proposition} \no $Proof:$ Uniqueness is already given in Proposition \ref{unique}. The existence
will be obtained via a fixed point argument. Actually, let us set
${\cal D}:={\cal H}^{2,1}$ endowed with the norm
$$||Y||_{\alpha}=E[\integ{0}{T}e^{\alpha
s}|Y_s|^2ds]^{\frac{1}{2}};\, \alpha>0.$$ Let $\Phi$ be the map from
$\cal D$ into itself defined by $\Phi(Y)=\tilde Y$ where $(\tilde
Y,\tilde Z,\tilde V,\tilde K^{\pm})$ is the solution of the
reflected BSDE associated with $(\xi,f(t,Y_t),L,U)$. Let $Y'$ be
another element of $\cal D$ and $\Phi(Y')=\tilde Y'$. Note again
that there is a lack of integrability for $(\tilde Z,\tilde V)$ and
$(\tilde Z',\tilde V')$, then we need to proceed by localisation. So
let us introduce the following sequence of stopping times:
$$\forall k\geq 1,\,\,\tau_k:=\mbox{inf}\{t\geq 0;
\integ{0}{t}(|Z_s|^2+|Z'_s|^2)ds
+\integ{0}{t}\int_E(|V_s(e)|^2+|V'_s(e)|^2)\lambda(de)ds\geq
k\}\wedge T.$$ As we discussed in Proposition \ref{unique}, the
sequence is non-decreasing, of stationary type and converges to $T$.
Applying It\^{o}'s formula to $e^{\alpha s}(\tilde Y_s-\tilde
Y'_s)^2$ on $[0,\tau_k]$, we will get: for any $t\leq T$,
\begin{equation}\label{eqcontr}
\begin{array}{ll}e^{\alpha (t\wedge\tau_k)}(\tilde
Y_{t\wedge\tau_k}-\tilde Y'_{t\wedge\tau_k})^2+\alpha
\integ{t\wedge\tau_k}{\tau_k}e^{\alpha s}(\tilde Y_s-\tilde
Y'_s)^2ds\\\leq
(M_{\tau_k}-M_{t\wedge\tau_k})+2\integ{t\wedge\tau_k}{\tau_k}e^{\alpha
s}(\tilde Y_{s-}-\tilde Y'_{s-})(d\tilde K_s^+-d\tilde K_s^{'+})
-2\integ{t\wedge\tau_k}{\tau_k}e^{\alpha s}(\tilde Y_{s-}-\tilde
Y'_{s-})(d\tilde K_s^--d\tilde K_s^{'-})\\+e^{\alpha \tau_k}(\tilde
Y_{\tau_k}-\tilde
Y'_{\tau_k})^2+2\integ{t\wedge\tau_k}{\tau_k}e^{\alpha s}(\tilde
Y_{s}-\tilde Y'_{s})(f(s,Y_s)-f(s,Y'_s))ds,
\end{array}\end{equation}
where $(M_{t\wedge\tau_k})_{t\leq T}$ is actually a martingale. But
taking into account Remark \ref{rmqimp}-$(ii)$, we deduce that:
$$\begin{array}{ll}\integ{t\wedge\tau_k}{\tau_k}e^{\alpha s}(\tilde
Y_{s-}-\tilde Y'_{s-})(d\tilde K_s^+-d\tilde
K_s^{'+})&=\integ{t\wedge\tau_k}{\tau_k}e^{\alpha s}(\tilde
Y_{s-}-S_{s-}+S_{s-}-\tilde Y'_{s-})(d\tilde K_s^+-d\tilde
K_s^{'+})\\{}&=\integ{t\wedge\tau_k}{\tau_k}e^{\alpha
s}(S_{s-}-\tilde Y'_{s-})d\tilde K_s^+-
\integ{t\wedge\tau_k}{\tau_k}e^{\alpha s}(\tilde
Y_{s-}-S_{s-})d\tilde K_s^{'+}\leq 0.\end{array}$$On the other hand,
since $(\tau_k)_{k\geq 1}$ is stationary, we have that $e^{\alpha
\tau_k}(\tilde Y_{\tau_k}-\tilde Y'_{\tau_k})^2\rightarrow 0$ when
$k\rightarrow \infty$. Therefore taking expectation in both hand
sides of (\ref{eqcontr}), using the inequality $|a.b|\leq
\epsilon^{-1}|a|^2 +\epsilon |b|^2 $ for any $\epsilon >0$ and
$a,b\in R^p$, and passing to the limit as $k\rightarrow \infty$, we
get:
$$
\begin{array}{ll}
(\alpha-\epsilon C_f)E[\integ{t}{T}e^{\alpha s}(\tilde Y_s-\tilde
Y'_s)^2ds]\leq \frac{C_f}{\epsilon}E[\integ{t}{T}e^{\alpha
s}(Y_s-Y'_s)^2ds].
\end{array}
$$
Choose $\alpha$ and $\epsilon$ appropriately, we can make that
$\Phi$ is a contraction on $\cal D$. Therefore it has a fixed point
$Y$ which belongs also to ${\cal S}^2$. Thus the proposition is
proved. $\Box$
\medskip

Additionally, we have also the following lemma related to local
integrability of the processes $Z,V$ and $K^\pm$.
\begin{lemma}\label{inte2} Assume [H] and let $(Y,Z,V,K^+,K^-)$ be the unique solution associated
with\\
  $(f(t,y),\xi,L,U)$. Let $(\gamma_k)_{k\geq 0}$ be a sequence of stopping
  times which satisfies $(i)$ and $(ii)$ of
  Proposition \ref{weakm}. Then for any $k\geq 0$, we have:
  $$E[\integ{0}{\gamma_{k}}ds\{|Z_s|^2+\int_E|V_s(e)|^2\lambda(de)\}+(K^+_{\gamma_k})^2+(K^-_{\gamma_k})^2]<\infty.$$
\end{lemma}
\no $Proof:$ Since the 5-uple $(Y,Z,V,K^+,K^-)$ is the solution of
the BSDE associated with $(f(t,y),\xi,L,U)$, then for any
$\gamma_k$, we have:
\begin{equation}\label{eqbsdeint}Y_{t\wedge\gamma_k}=Y_{\gamma_k}+\integ{t\wedge\gamma_k}{\gamma_k}f(s,Y_s)ds
+\integ{t\wedge\gamma_k}{\gamma_k}d(K_s^{+}-K_s^{-})
-\integ{t\wedge\gamma_k}{\gamma_k}Z_sdB_s-\integ{t\wedge\gamma_k}{\gamma_k}
\int_EV_s(e)\widetilde\mu(ds,de).\end{equation} On the other hand,
on $[0,\gamma_k]$, Mokobodzki's condition
 \bf{[Mk]} is satisfied. Therefore the BSDE associated with
$(f(t,y)1_{[t\leq\gamma_k]},Y_{\gamma_k},L_{t\wedge\gamma_k},U_{t\wedge\gamma_k})$
has a solution (see e.g. \cite{ho}) which we denote \\
$(Y^k,Z^k,V^k,K^{k,+},K^{k,-})$. Then it holds true that for any
$t\leq T$:
$$Y_{t\wedge\gamma_k}^k=Y_{\gamma_k}+\integ{t\wedge\gamma_k}{\gamma_k}f(s,Y_s^k)ds
+\integ{t\wedge\gamma_k}{\gamma_k}d(K_s^{k,+}-K_s^{k,-})
-\integ{t\wedge\gamma_k}{\gamma_k}Z_s^kdB_s-\integ{t\wedge\gamma_k}{\gamma_k}
\int_EV_s^k(e)\widetilde\mu(ds,de).$$ Moreover we have the following
integrability property:
\begin{equation}\label{inte3}
E[\integ{0}{\gamma_{k}}ds\{|Z_s^k|^2+\int_E|V_s^k(e)|^2\lambda(de)\}+(K^{k,+}_{\gamma_k})^2+(K^{k,-}_{\gamma_k})^2]<\infty.
\end{equation}
But uniqueness of the solution of the BSDE (\ref{eqbsdeint}) implies
that:
$$Y_{t\wedge\gamma_k}=Y_{t\wedge\gamma_k}^k,Z_{t\wedge\gamma_k}=Z_{t\wedge\gamma_k}^k,
V_{t\wedge\gamma_k}=V_{t\wedge\gamma_k}^k,
K_{t\wedge\gamma_k}^{+}-K_{t\wedge\gamma_k}^{-}=K_{t\wedge\gamma_k}^{k,+}-K_{t\wedge\gamma_k}^{k,-}.$$
Therefore, the desired result follows from (\ref{inte3}). $\Box$
\bigskip

We are now ready to establish the main result of this paper. The
proof is basically the same as the one given in (\cite{hh3}, Theorem
4.2, Step 2) even if in this latter paper the obstacles have only
inaccessible jumps, therefore it is omitted.
\begin{theorem}
\label{thmprince} Under \bf{[H]}, the BSDE (\ref{bsde}) with jumps
and two reflecting discontinuous barriers associated with $(f,\xi,
L,U)$ has a unique solution, i.e., there exits a unique 5-uple
$(Y,Z,V,K^+,K^-)$ which satisfies the BSDE (\ref{bsde}). $\Box$
\end{theorem}

\begin{remark}
Under $[H]$ we have also the uniqueness of the increasing processes.
Actually if \\ $(Y,Z,V,K^\pm)$ and $(Y',Z',V',K'^\pm)$ are two
solutions of the BSDE associated with $(f(t,y,z,v),\xi,L,U)$ then we
have also $K^+=K'^+$ and $K^-=K'^-$ (see e.g. \cite{hh3} for the
proof of this claim). $\Box$
\end{remark}

We now deal with an application of these types of BSDEs in zero-sum
mixed game problems.
\section{Application in zero-sum mixed differential-integral game problem}
We are going now to study the link between mixed zero-sum stochastic
differential game and the reflected BSDE studied in the previous
section. First let us briefly describe the setting of the problem of
zero-sum game we consider.

Let $x_0\in R^d$ and let $x=(x_t)_{t\leq T}$ be the solution of the
following standard differential equation:
$$x_t=x_0+\integ{0}{t}\sigma(s,x_s)dB_s+\integ{0}{t}\int_E\gamma(s,e,x_{s-})\widetilde\mu(ds,de)$$
where the mapping $\sigma$:
$(t,x)\in[0,T]\times{R^d}\mapsto\sigma(t,x)\in R^d$ and $\gamma$:
$(t,e,x)\in [0,T]\times E\times R^d \mapsto\gamma(t,e,x)\in R^d$
satisfy
the following assumptions: \\
(i): there exists a constant $C_1$ such that
$$\forall (t,x), tr
(\sigma\sigma^\ast(t,x))+\integ{E}{}\gamma(t,e,x)^2\lambda(de)\leq
C_1(1+|x|^2);$$ (ii): there exists a constant $C_2$ such that
$$\forall (t,x),tr[(\sigma(t,x)-\sigma(t,y))(\sigma^\ast(t,x)-\sigma^\ast(t,y))]
+\integ{E}{}|\gamma(t,e,x)-\gamma(t,e,y)|^2\lambda(de)\leq
C_2|x-y|^2;$$ (iii): $\forall (t,x)\in [0,T]\times R^d$, the matrix
$\sigma(t,x)$ is invertible and $\sigma^{-1}(t, x)$ is bounded.
\medskip

According to Theorem 1.19 in \cite{sulok}, the process $(x_t)_{t\leq
T}$ exists and is unique thanks to the assumptions (i)-(ii) on the
functions $\sigma$ and $\gamma$. $\Box$
\medskip

Let $A\,\,(resp. \,B)$ be a compact metric space and ${\cal
U}\,\,(resp.\, \cal V)$ be the space of $\cal P$-measurable
processes $u=(u_t)_{t\leq T}\,\,(resp.\,v=(v_t)_{t\leq T})$ with
values in $A\,\,(resp. \,B)$. Let $f$ be a function from
$[0,T]\times{R^d}\times A\times B$ into $R^d$ which is ${\cal
B}([0,T]\times R^d\times A\times B)$-measurable and which satisfies:

(3-a): $f(t,x,u,v)$ is bounded for any $t,x,u \mbox{ and }v$;

(3-b): for any $(t,x)\in[0,T]\times{R^d}$, the mapping $(u,v)\mapsto
f(t,x,u,v)$ is continuous .$\Box$
\medskip

Now for $(u,v)=(u_t,v_t)_{t\leq T}\in {\cal U\times V}$, let
$L^{u,v}:=(L_t^{u,v})_{t\leq T}$ be the positive local martingale
solution of:
$$dL_t^{u,v}=L_{t-}^{u,v}\{\sigma^{-1}(t,x_t)f(t,x_t,u_t,v_t)dB_t
 +\int_E\beta(t,e,x_{t-},u_t,v_t)\widetilde\mu(dt,de)\}\mbox{ and }L_0^{u,v}=1$$
where for any $t,x,e,u,v$ we have $-1<\beta(t,x,e,u,v)$ and
$|\beta(t,x,e,u,v)|\leq c_0 (1\wedge |e|)$ where $c_0$ is a
constant. Then the measure $P^{u,v}$ defined by:
$$\frac{dP^{u,v}}{dP}{|}_{{\cal F}_T}=L_T^{u,v}$$
is actually a probability (\cite{bl}, Corollary 5.1, pp.244)
equivalent to $P$. Moreover, under the new probability $P^{u,v}$,
$\mu(dt,de)$ remains a random measure, whose compensator is
\\$\bar\nu(dt,de)=(1+\beta(t,e,x_{t-},u_t,v_t))\lambda(de)dt$, i.e.
$\widetilde\mu^{u,v}([0,t]\times A):=(\mu-\bar\nu)([0,t]\times
A)_{t\leq T}\,$ is a martingale for any $A\in\cal E$ satisfying
$\lambda(A)<\infty$, and
$B_t^{u,v}=B_t-\int_{0}^{t}\sigma^{-1}(s,x_s)f(s,x_s,u_s,v_s)ds$ is
a Brownian motion and $(x_t)_{t\leq T}$ satisfies:
$$\begin{array}{ll}
 x_t&=x_0+\integ{0}{t}f(s,x_s,u_s,v_s)ds+\integ{0}{t}\sigma(s,x_s)dB_s^{u,v}+
            \integ{0}{t}\!\!\!\int_E\gamma(s,e,x_{s-})\widetilde\mu^{u,v}(ds,de)\\
    &\qquad \qquad+\integ{0}{t}\!\!\!\int_E\gamma(s,e,x_{s-})\beta(s,e,x_{s-},u_s,v_s)\lambda(de)ds
\end{array}$$
It means that $(x_t)_{t\leq T}$ is a weak solution for this
stochastic differential equation and it stands for the evolution of
a system when controlled.\medskip

As we know, in mixed game problems, on a system intervene two agents
$c_1$ and $c_2$ who act with admissible controls $u$ and $v$
respectively which belong to $\cal U$ and $\cal V$ respectively.
Moreover, they can make the decision to stop controlling at $\tau$
for $c_1$ and $\sigma$ for $c_2$, where $\tau$ and $\sigma$ are two
stopping times. Therefore a strategy for $c_1$ (resp. $c_2$) is a
pair $(u,\tau)$ (resp. $(v,\sigma)$) and the system is actually
stopped at $\tau\wedge\sigma$. Meanwhile, the interventions of the
agents will generate a payoff which is a cost for $c_1$ and a reward
for $c_2$ whose expression is given by:
$$J(u,\tau;v,\sigma)=E^{u,v}[\integ{0}{\tau\wedge\sigma}h(s,x_s,u_s,v_s)ds+U_{\tau}1_{[\tau<\sigma]}
+L_{\sigma}1_{[\sigma<\tau<T]}+\xi1_{[\tau=\sigma=T]}],$$ where:

(1): $h:[0,T]\times R^d\times A\times B \mapsto R^+$ is ${\cal
P}\bigotimes{{\cal B}(A\times B)}$-measurable function which stands
for the instantaneous payoff between the two agents. In addition,
the mapping is continuous w.r.t. $(u,v)$ and there exists a constant
$C_h$ such that for any $(t,x,u,v)$, $|h(t,x,u,v)|\leq C_h(1+|x|)$;

(2): the stopping payoffs $U=(U_t)_{t\leq T}$ and $L=(L_t)_{t\leq
T}$ are processes of ${\cal S}^2$ and satisfy assumption [H],
$i.e.$, $L_t<U_t$ and $L_{t-}<U_{t-}$ $\forall$ $t\leq T$;

(3): $\xi$ is a ${\cal F}_T$-measurable random variable such that
$E[\xi^2]<\infty$ and $L_T\leq\xi\leq U_T$.
\begin{remark}
Here we assume that $L$ and $U$ are strictly separated in order to
infer the existence of a global solution of a RBSDE associated with
$\xi$, $L$, $U$ and an appropriate generator which we will precise
later. $\Box$
\end{remark}

In this zero-sum game problem we aim at showing that the value of
the game exists, i.e., it holds true that:
\begin{equation}\label{valeur}
essinf_{(u,\tau)}esssup_{(v,\sigma)}J(u,\tau;v,\sigma)=
esssup_{(v,\sigma)}essinf_{(u,\tau)}J(u,\tau;v,\sigma).\end{equation}

In \cite{[H2],[HL1]}, the authors deal the mixed zero-sum
differential game when the process $(x_t)_{t\leq T}$ has no jump
part,  the information comes only from a Brownian motion and the
stopping payoffs are continuous. Actually, using results on two
barrier reflected BSDEs they proved that the zero-sum game has a
value and a saddle-point also. The value is expressed by means of
the solution of the BSDE with two reflecting barriers. In this work,
and for our general setting, we will be just able to show that the
value of the mixed zero-sum differential game exits. However we are
not able to infer the existence of a saddle-point because, and this
is the main reason for that, the payoff processes have predictable
jumps.
\medskip

So let us define the Hamilton function associated with this game
problem as following: \\ $\forall (t,x,z,r,u,v)\in [0,T]\times
R^d\times R^d\times L_R^2(E,d\lambda)\times A\times B$,
$$H(t,x,z,r,u,v):=z\sigma^{-1}(t,x)f(t,x,u,v)
+h(t,x,u,v)+\int_E r(e)\beta(t,e,x,u,v)\lambda(de)$$ Next assume
that Isaacs'condition, which plays an important role in zero-sum
mixed game problems, is fulfilled, $i.e.$, for any
$(t,x,z,r)\in[0,T]\times R^d\times R^d\times L_R^2(E,d\lambda)$,
$$[A]: \textrm{inf}_{u\in A}\textrm{sup}_{v\in B}H(t,x,z,r,u,v)=
\textrm{sup}_{v\in B}\textrm{inf}_{u\in A}H(t,x,z,r,u,v).$$

Under [A], through the assumptions above and Benes' selection
theorem, the following result holds true (see e.g. \cite{elkham}).
\begin{proposition}
There exist two measurable functions $u^\ast(t,x,z,r)$ and
$v^\ast(t,x,z,r)$ from\\ $[0,T]\times R^d\times R^d\times
L_R^2(E,d\lambda)$ into $A$ and $B$ respectively, such that:

(i) the pair $(u^\ast,v^\ast)(t,x,z,r)$ is a saddle-point for the
function $H$, i.e., for any $u,v$ we have:
$$H(t,x,z,r,u^\ast(t,x,z,r),v)\leq H(t,x,z,r,
(u^\ast,v^\ast)(t,x,z,r))\leq H(t,x,z,r,u,v^\ast(t,x,z,r)).$$

(ii) the function $(z,r)\mapsto H(t,x,z,r,(u^\ast,v^\ast)(t,x,z,r))$
is uniformly Lipschitz.
\end{proposition}
\medskip

Now let us set
$H^\ast(t,x_t(\omega),z,r)=H(t,x_t(\omega),z,r,(u^\ast,v^\ast)(t,x_t(\omega),z,r))$
and let $(Y_t,Z_t,R_t,K_t^{\pm})$ be the global solution associated
with $(H^\ast,\xi,L,U)$, which exists according to Theorem
\ref{thmprince}. Therefore we have: $\forall t\in[0,T] $,
\begin{equation}\label{Hstar}
 \left\{\!\!\begin{array}{ll}
   (i)Y_t=\xi+\integ{t}{T}\!\!H^\ast(s,x_s,Z_s,R_s)ds+(K^+_T-K^+_t)
   -(K^-_T-K^-_t)-\integ{t}{T}Z_sdB_s-\integ{t}{T}\!\!\!\int_ER_{s}(e)\widetilde\mu(ds,de)   \\
      (ii) \,\,L_t\leq Y_t\leq U_t,\,  \integ{0}{T} (U_s-Y_s)dK_s^{
c,-}=\integ{0}{T} (Y_s-L_s)dK_s^{c,+}=0 \mbox{ where }K^{
c,\pm} \mbox{ is the }\\
\qquad\mbox{continuous part of } K^{\pm}\, (K^{
c,+}_0=0);\\
(iii)\,\, K^{d,\pm}, \mbox{ the purely discontinuous part of
}K^{\pm}\mbox{ is predictable and verifies }\\
\qquad\qquad \qquad K^{d,+}_t=\sum_{0<s\leq t}(L_{s-}-Y_s)^+ \mbox{
and }K^{
d,-}_t=\sum_{0<s\leq t}(Y_s-U_{s-})^+;\\
(iv) \integ{0}{T}|Z_s|^2ds+
\integ{0}{T}\!\!\!\int_E|R_s(e)|^2\lambda(de)ds<\infty, P-a.s.
       \end{array}
 \right.
\end{equation}

The following is the main result of this part:
\begin{theorem}\label{mixedgame} We have:
$$Y_0=\mbox{esssup}_{\sigma\in{{\T}_0},v\in{\cal
V}}\mbox{essinf}_{\tau\in{{\T}_0},u\in{\cal U}}J(u,\tau;v,\sigma)=
\mbox{essinf} _{\tau\in{{\T}_0},u\in{\cal
U}}\mbox{esssup}_{\sigma\in{{\T}_0},v\in{\cal
V}}J(u,\tau;v,\sigma)$$i.e. $Y_0$ is the value of the zero-sum mixed
differential game.
\end{theorem}
$Proof$: First note that $Y_0$ is a constant since ${\cal F}_0$
contains only the $P$-null sets of ${\cal F}$. Now, for any fixed
$(u,v)\in{\cal U}\times{\cal V}$, let $(Y^{u,v},\bar Z, \bar R, \bar
K^\pm)$ be the solution of the following reflected BSDE:
\begin{equation}\label{H}
 \left\{\begin{array}{ll}
   (i)\,\,Y_t^{u,v}=\xi+\int_{t}^{T}H(s,x_s,\bar Z_s,\bar R_s,u_s,v_s)ds+(\bar K^+_T-\bar K^+_t)
   -(\bar K^-_T-\bar K^-_t)\\\qquad\qquad\qquad-\int_{t}^{T}\bar Z_sdB_s-\int_{t}^{T}\!\!\!\int_E\bar R_s(e)\widetilde\mu(ds,de);
  \\ (ii)\,\,  L_t\leq Y_t^{u,v}\leq U_t\, and \int_{0}^{T} (U_s-Y_s^{u,v})d\bar K_s^{
c,-}=\int_{0}^{T} (Y_s^{u,v}-L_s)d\bar K_s^{c,+}=0 \mbox{ where
}\bar K^{ c,\pm} \\\qquad\mbox{ is the continuous part of } \bar
K^{\pm}\,\, (\bar K^{
c,+}_0=0);\\
(iii)\,\,\bar K^{d,+}_t=\sum_{0<s\leq t}(L_{s-}-Y_s^{u,v})^+ \mbox{
and }\bar
K^{ d,-}_t=\sum_{0<s\leq t}(Y_s^{u,v}-U_{s-})^+;\\
(iv)\,\,\int_{0}^{T}|\bar Z_s|^2ds+ \int_{0}^{T}\!\!\!\int_E|\bar
R_s(e)|^2\lambda(de)ds<\infty, P-a.s.
       \end{array}
 \right.
\end{equation}
Even in our setting where there are general jumps in the equation,
making a change of probability and arguing as in \cite{tmx}, we
obtain that $Y_t^{u,v}$ is the value function of the Dynkin game,
$i.e.$,
$$Y_t^{u,v}=\mbox{esssup}_{\sigma\in{{\T}_t}}\mbox{essinf}_{\tau\in{{\T}_t}}J_t(u,\tau;v,\sigma)=
\mbox{essinf}
_{\tau\in{{\T}_t}}\mbox{esssup}_{\sigma\in{{\T}_t}}J_t(u,\tau;v,\sigma),$$
where
\begin{equation}\label{cout}J_t(u,\tau;v,\sigma)=E^{u,v}[\int_{t}^{\tau\wedge\sigma}h(s,x_s,u_s,v_s)ds+U_{\tau}1_{[\tau<\sigma]}
+L_{\sigma}1_{[\sigma\leq\tau<T]}+\xi1_{[\tau=\sigma=T]}|{\cal
F}_t].\end{equation}

Let us now prove that:
\begin{equation}\label{val}Y_t=\mbox{esssup}_{v\in{\cal V}}\mbox{essinf}_{u\in\cal
U}Y_t^{u,v}=\mbox{essinf}_{u\in\cal U}\mbox{esssup}_{v\in{\cal
V}}Y_t^{u,v}.\end{equation} However since $\mbox{esssup}_{v\in{\cal
V}}\mbox{essinf}_{u\in{\cal U}}Y_t^{u,v}\leq\mbox{essinf}_{u\in{\cal
U}}\mbox{esssup}_{v\in{\cal V}}Y_t^{u,v}$, we just need to prove
that:
$$\mbox{essinf}_{u\in{\cal U}}\mbox{esssup}_{v\in\cal
V}Y_t^{u,v}\leq Y_t\leq\mbox{esssup}_{v\in\cal
V}\mbox{essinf}_{u\in{\cal U}}Y_t^{u,v}$$ where $Y_t$ is the
solution of (\ref{Hstar}).

First note that the processes
$(u^\ast_t=u^\ast(t,x_t,Z_t,R_t))_{t\leq T}$ and
$(v^\ast_t=v^\ast(t,x_t,Z_t,R_t))_{t\leq T}$ are admissible
controls. Now let $(u_t)_{t\leq T}$ be an arbitrary admissible
control. The generator \\$H(t,x_t,z,r,u_t,v^\ast(t,x_t,Z_t,R_t))$ is
uniformly Lipschitz $w.r.t.$ $(z,r)$. Therefore thanks to Theorem
\ref{thmprince} there exists a process $Y^{u,v^\ast}$ such that for
any $t\leq T$:
$$Y_t^{u,v\ast}=\xi+\integ{t}{T}H(s,\tilde Z_s,\tilde R_s,u_s,v_s^{\ast})ds+(\tilde K^+_T-\tilde K^+_t)
   -(\tilde K^-_T-\tilde K^-_t)-\integ{t}{T}\tilde Z_sdB_s-\integ{t}{T}\!\!\!\int_E\tilde R_s(e)\widetilde\mu(ds,de).$$
Let us define a new probability $P^{u,v^\ast}$ by
$$\frac{dP^{u,v^\ast}}{dP}{|}_{{\cal F}_T}=L_T^{u,v^\ast}.$$
Using It\^{o}-Meyer's formula (\cite{protter}, pp.221) for
${(Y-Y^{u,v^{\ast}})^+}^2$ and taking into account that:
$$\begin{array}{ll}
          H^\ast(s,x_s,Z_s,R_s)&-H(s,x_s,\tilde Z_s,\tilde R_s,u_s,v^\ast_s)=H^\ast(s,x_s,Z_s,R_s)-H(s,x_s,Z_s,R_s,u_s,v^\ast_s)\\
         &\qquad\qquad +H(s,x_s,Z_s,R_s,u_s,v^\ast_s)-H(s,x_s,\tilde
         Z_s,\tilde R_s,u_s,v^\ast_s)\\
         &=H^\ast(s,x_s,Z_s,R_s)-H(s,x_s,Z_s,R_s,u_s,v^\ast_s)+(Z_s-\tilde
         Z_s)\sigma^{-1}(s,X_s)f(s,X_s,u_s,v^\ast_s)\\&\quad\qquad +\integ{E}{}(R_{s-}(e)-\tilde
         R_{s-}(e))\beta(s,e,X_s,u_s,v^\ast_s)\lambda(de),
\end{array}$$
we obtain: $\forall t\in[0,T]$:
 $$\begin{array}{ll}
     {(Y_t-Y_t^{u,v^{\ast}})^+}^2&\leq2\int_{t}^{T}(Y_s-Y_s^{u,v^{\ast}})^+
     (H^\ast(s,x_s,Z_s,R_s)-H(s,x_s,Z_s,R_s,u_s,v^\ast_s))ds\\
     &+2\int_{t}^{T}(Z_s-\tilde
     Z_s)dB_s^{u,v^\ast}+2\int_{t}^{T}\int_E(R_s(e)-\tilde
     R_s(e))\widetilde\mu^{u,v^\ast}(ds,de),
\end{array}$$
where under the new probability $P^{u,v^\ast}$, the process
$B^{u,v^\ast}$ is a Brownian motion and $\mu^{u,v^\ast}(ds,de)$ is a
martingale measure. Now since
$H^\ast(s,x_s,Z_s,R_s)-H(s,x_s,Z_s,R_s,u_s,v^\ast_s)\leq0$, after
localization, taking expectation under $P^{u,v^\ast}$ and then the
limit, we obtain $P^{u,v^\ast}-a.s.$, $Y_t\leq Y_t^{u,v^\ast}.$
Therefore $P-a.s.$ for any $t\leq T$, $Y_t\leq Y_t^{u,v^\ast}$ since
the two probabilities are equivalent. In the same way we can show
that $Y_t^{u^\ast,v}\leq Y_t$, $P-a.s.$ for any $t\leq T$ and any
admissible control $(v_t)_{t\leq T}$. Therefore for any $t\leq T$ we
have:
$$Y_t^{u^\ast,v}\leq Y_t\leq Y_t^{u,v^\ast}$$ and then
$$\mbox{essinf}_{u\in{\cal U}}\mbox{esssup}_{v\in\cal
V}Y_t^{u,v}\leq Y_t\leq \mbox{esssup}_{v\in\cal
V}\mbox{essinf}_{u\in\cal U}Y_t^{u,v}$$ which ends the proof of
(\ref{val}).

We now focus on the main claim. So let us prove that:
\begin{equation}\label{commute}\mbox{essinf}_{u\in\cal U}Y_t^{u,v}
=\mbox{esssup}_{\sigma\in{{\T}_t}}\mbox{essinf}_{\tau\in{{\T}_t}}\mbox{essinf}_{u\in{\cal
U}}J_t(u,\tau;v,\sigma)\end{equation} $i.e.$ we can commute the
control and the stopping times. So for any $u,v$ and $\sigma$,
$\tau$ let $(J_t,Z_t,r_t)_{t\leq\tau\wedge\sigma}$ be the solution
of the following standard BSDE:
$$J_t=\bar\xi+\integ{t}{\tau\wedge\sigma}H(s,z_s,r_s,u_s,v_s)ds-
\integ{t}{\tau\wedge\sigma}z_sdB_s-\integ{t}{\tau\wedge\sigma}\int_E
r_s(e)\widetilde\mu(ds,de)$$ where
$\bar\xi=U_{\tau}1_{[\tau<\sigma]}
+L_{\sigma}1_{[\sigma\leq\tau<T]}+\xi1_{[\tau=\sigma=T]}$. This
solution exists thanks to a result by Tang $\&$ Li \cite{tangli}.
Therefore $P-a.s.$, for any $t\leq \tau \wedge \sigma$, we have
$J_t= J_t(u,\tau;v,\sigma)$.

We can now argue as in \cite{[EKPQ]}, Proposition 3.1, to obtain
that:
$$\mbox{essinf}_{u\in{\cal U}}J_t(u,\tau;v,\sigma)=\bar\xi+\integ{t}{\tau\wedge\sigma}
\mbox{essinf}_{u\in{\cal U}}H(s,z_s,r_s,u_s,v_s)ds-
\integ{t}{\tau\wedge\sigma}z_sdB_s-\integ{t}{\tau\wedge\sigma}\!\!\int_E
r_s(e)\widetilde\mu(ds,de).$$ Actually this is possible since we can
use comparison of solutions of those BSDEs thanks to the properties
satisfied by the mapping $\beta$ and especially the fact that $\beta
>-1$. Therefore the process\\
$(\mbox{esssup}_{\sigma\in{{\T}_t}}\mbox{essinf}_{\tau\in{{\T}_t}}\mbox{essinf}_{u\in{\cal
U}}J_t(u,\tau;v,\sigma))_{t\leq T}$ is the value function of the
corresponding Dynkin game, $i.e.$ the solution of the RBSDE
associated with $(\mbox{essinf}_{u\in{\cal
U}}H(t,z,r,u,v),\xi,L,U)$.
\medskip

On the other hand,  once more using comparison of solutions of BSDEs
with two reflecting barriers we obtain that the process
$(\mbox{essinf}_{u\in\cal U}Y_t^{u,v})_{t\leq T}$ is the solution
(with the other components) of the RBSDE associated
$(\mbox{essinf}_{u\in{\cal U}}H(t,z,r,u,v),\xi,L,U)$. Now by
uniqueness we obtain: for any $t\leq T$,
$$
\mbox{essinf}_{u\in\cal
U}Y_t^{u,v}=\mbox{esssup}_{\sigma\in{{\T}_t}}\mbox{essinf}_{\tau\in{{\T}_t}}\mbox{essinf}_{u\in{\cal
U}}J_t(u,\tau;v,\sigma).$$ It follows that: $\forall t\leq T$,
$$\begin{array}{ll}
Y_t&=\mbox{essup}_{v\in\cal V}\mbox{essinf}_{u\in\cal
U}Y_t^{u,v}=\mbox{essup}_{v\in\cal
V}\mbox{esssup}_{\sigma\in{{\T}_t}}\mbox{essinf}_{\tau\in{{\T}_t}}\mbox{essinf}_{u\in{\cal
U}}J_t(u,\tau;v,\sigma)\\{}&=\mbox{esssup}_{\sigma\in{{\T}_0},v\in{\cal
V}}\mbox{essinf}_{\tau\in{{\T}_0},u\in{\cal
U}}J_t(u,\tau;v,\sigma).\end{array}$$ In the same way we can show
that:
$$
\mbox{esssup}_{v\in\cal
V}Y_t^{u,v}=\mbox{essinf}_{\tau\in{{\T}_t}\mbox{esssup}_{\sigma\in{{\T}_t}}}\mbox{esssup}_{v\in{\cal
V}}J_t(u,\tau;v,\sigma)$$which implies that:
$$\begin{array}{ll}
Y_t=\mbox{essinf}_{\tau\in{{\T}_0},u\in{\cal
U}}\mbox{esssup}_{\sigma\in{{\T}_0},v\in{\cal
V}}J_t(u,\tau;v,\sigma),\,\,t\leq T.\end{array}$$ Thus the proof of
the claim is complete. $\Box$
\begin{center}\bf{ Appendix}
\end{center}
\begin{lemma}\label{l1}: If $(U^n)_{n\geq 0}$ is a non-decreasing sequence of
progressively measurable $rcll$ $\R$-valued processes of class [D]
which converges pointwisely to $U$ another progressively measurable
$rcll$ $\R$-valued process of class [D], then $P$-a.s., $\forall
t\leq T$, \textsl{SN}$(U^n)_t\nearrow\textsl{SN}(U)_t$, where
\textsl{SN} is the Snell envelope operator.
\end{lemma}
The proof of this result has been given in several works (see e.g.
Appendix in \cite{ho}) and then we omit it. $\Box$

\small{
}
\end{document}